\documentclass[aos,MSNbibl,nameyear,dvips]{arximspdf}
\usepackage{accents}
\usepackage{graphicx}
\doi{10.1214/14-AOS1218} 
\volume{42}
\issue{4}
\pubyear{2014}
\firstpage{1511}
\lastpage{1545}

\makeatletter
\newproclaim{example}{Example}
\newproclaim{remark}{Remark}

\def\VHF{\approx^{\hspace*{-8.6pt}\fontsize{11}{11}\selectfont{\mbox{$\sim$}}}}
\def\cal{\mathcal}

\def\tilde{\widetilde}

\newtheorem{theorem}{Theorem}
\newtheorem{lemma}[theorem]{Lemma}
\newtheorem{corollary}[theorem]{Corollary}

\newcommand\cK{{\cal K}}
\newcommand\cM{{\cal M}}

\def\R{\mathbb{R}}
\newcommand{\Hausdorff}{\operatorname{Hausdorff}}

\newcommand{\dt}[1]{\accentset{\mbox{\textbf{\large .}}}{#1}}
\makeatother

\begin{document}
\begin{frontmatter}

\title{Nonparametric ridge estimation}
\runtitle{Ridge estimation}

\begin{aug}
\author[a]{\fnms{Christopher R.}~\snm{Genovese}\thanksref{t1}\ead[label=e1]{genovese@stat.cmu.edu}},
\author[b]{\fnms{Marco}~\snm{Perone-Pacifico}\thanksref{t2}\ead[label=e2]{marco.peronepacifico@uniroma1.it}},
\author[c]{\fnms{Isabella}~\snm{Verdinelli}\thanksref{t2}\ead[label=e3]{isabella@stat.cmu.edu}}
\and
\author[a]{\fnms{Larry}~\snm{Wasserman}\corref{}\thanksref{t4}\ead[label=e4]{larry@stat.cmu.edu}}
\thankstext{t1}{Supported by NSF Grant DMS-08-06009.}
\thankstext{t2}{Supported by Italian National Research Grant PRIN 2008.}
\thankstext{t4}{Supported by NSF Grant DMS-08-06009, Air Force Grant
FA95500910373.}
\runauthor{Genovese, Perone-Pacifico, Verdinelli and Wasserman}
\affiliation{Carnegie Mellon University, Sapienza University of Rome,
Carnegie Mellon University and Sapienza University of Rome,
and Carnegie Mellon University}
\address[a]{C.~R. Genovese\\
L. Wasserman\\
Department of Statistics\\
Carnegie Mellon University\\
Pittsburgh, Pennsylvania 15213\\
USA\\
\printead{e1}\\
\phantom{E-mail:\ }\printead*{e4}}

\address[b]{M. Perone-Pacifico\\
Department of Statistical Sciences\\
Sapienza University of Rome\\
Rome\\
Italy\\
\printead{e2}}

\address[c]{I. Verdinelli\\
Department of Statistics\\
Carnegie Mellon University\\
Pittsburgh, Pennsylvania 15213\\
USA\\
and\\
Department of Statistical Sciences\\
Sapienza University of Rome\\
Rome\\
Italy\\
\printead{e3}}
\end{aug}

\received{\smonth{6} \syear{2013}}
\revised{\smonth{3} \syear{2014}}

%
\begin{abstract}
We study the problem of estimating the ridges of a density function.
Ridge estimation is an extension of mode finding and
is useful for understanding the structure of a density.
It can also be used to find
hidden structure in point cloud data.
We show that, under mild regularity conditions,
the ridges of the kernel density estimator
consistently estimate the ridges
of the true density.
When the data are noisy measurements of a manifold,
we show that the ridges are close and topologically similar to the
hidden manifold.
To find the estimated ridges in practice, we adapt
the modified mean-shift algorithm proposed by
Ozertem and Erdogmus [\textit{J. Mach. Learn. Res.}
\textbf{12} (2011) 1249--1286].
Some numerical experiments verify that the algorithm
is accurate.
\end{abstract}

%
\begin{keyword}[class=AMS]
\kwd[Primary ]{62G05}
\kwd{62G20}
\kwd[; secondary ]{62H12}
\end{keyword}

\begin{keyword}
\kwd{Ridges}
\kwd{density estimation}
\kwd{manifold learning}
\end{keyword}
\end{frontmatter}

\section{Introduction}\label{sec1}

Multivariate data in many problems exhibit intrinsic lower dimensional
structure.
The existence of such structure is of great interest
for dimension reduction, clustering and improved statistical inference,
and the question of how to identify and characterize
this structure is the focus of active research.
A commonly used representation for low-dimensional structure is a
smooth manifold.
Unfortunately, estimating manifolds can be difficult even under mild
assumptions.
For instance, the rate of convergence for estimating a manifold with
bounded curvature
blurred by homogeneous Gaussian noise,
is \emph{logarithmic} [\citet{genovese2012geometry}],
meaning that an exponential amount of data are needed to attain a
specified level of accuracy.
In this paper, we offer a way to circumvent this problem.
We define an object, which we call a \emph{hyper-ridge} set
that can be used to approximate the low-dimensional structure in a data set.
We show that the hyper-ridge set captures the essential features of the
underlying low-dimensional structure
while being estimable from data at a polynomial rate.

\begin{figure}

\includegraphics{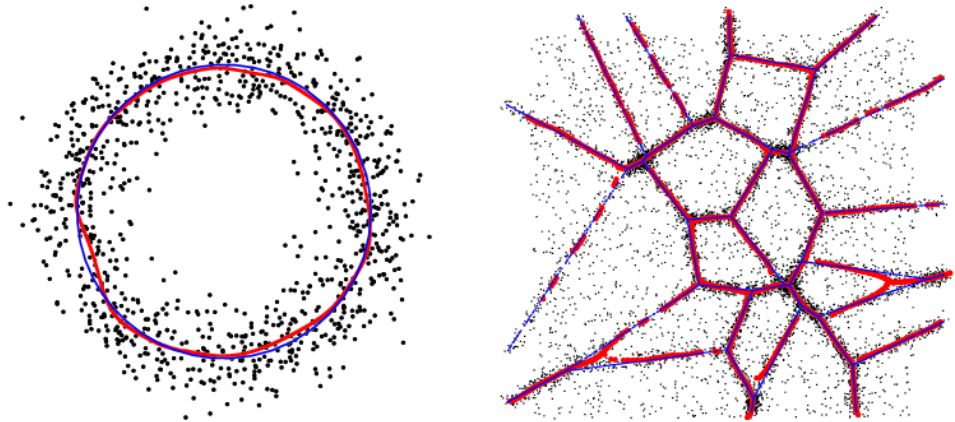}

\caption{Synthetic data showing lower dimensional structure.
The left plot is an example of the hidden manifold case.
The right plot is an example of a hidden set
consisting of intersecting manifolds.}
\label{fig::struct1}
\end{figure}

Let $X_1, \ldots, X_n$ be a sample from a probability density $p$
defined on an open subset of $D$-dimensional Euclidean space
and let $\hat p$ be an estimate of the density.
We will define hyper-ridge sets (called ridges for short)
for both $p$ and $\hat p$, which we denote by $R$ and $\hat R$.
We consider two cases that make different assumptions about $p$.
In the \emph{hidden manifold} case
(see Figure~\ref{fig::struct1}),
we assume that the density $p$ is derived by sampling
from a $d < D$ dimensional manifold $M$ and adding $D$-dimensional noise.
In the \emph{density ridge} case,
we look for ridges of a density without assuming any hidden manifold,
simply as a way of finding structure in a point cloud, much like clustering.
The goal in both cases is to estimate the hyper-ridge set.
Although in the former case, we would ideally like to estimate $M$,
this is not always
feasible for reasonable sample sizes, so we use the ridge $R$ as a
\emph
{surrogate}
for $M$.
We focus on estimating ridges from point cloud data;
we do not consider image data in this paper.

A formal definition of a ridge
is given in Section~\ref{sec::conditions}.
Let $1 \leq d < D$ be fixed.
Loosely speaking, we define a $d$-dimensional hyper-ridge set of a
density $p$
to be the points where the Hessian of $p$ has $D-d$ strongly negative
eigenvalues
and where the projection of the gradient on that subspace
is zero.
Put another way,
the ridge is a local maximizer of the density
when moving in the normal direction defined
by the Hessian.

Yet another way to think about ridges is by analogy with modes.
We can define a mode to be a point where the gradient is 0
and the second derivative is negative,
that is, the eigenvalues of the Hessian are negative.
The Hessian defines a $(D-d)$-dimensional normal space
(corresponding to the $D-d$ smallest eigenvalues)
and a $d$ dimensional tangent space.
A ridge point has a projected gradient
(the gradient in the direction of the normal)
that is 0
and eigenvalues in the normal space that are negative.
Modes are simply $0$ dimensional ridges.
\begin{figure}

\includegraphics{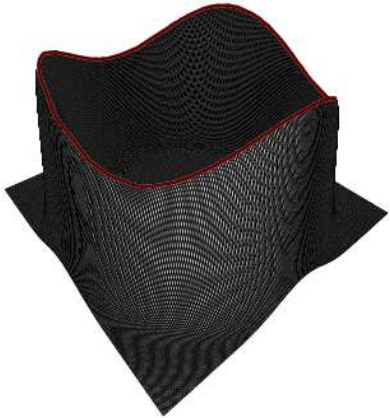}

\caption{An example of a one dimensional ridge defined by a
two-dimensional density $p$.
The ridge~$R$ is a circle on the plane.
The solid curve is the ridge, lifted onto $p$, that is,
$\{(x,p(x))\dvtx x\in R\}$.}
\label{fig::Pedagogical}
\end{figure}

\begin{figure}

\includegraphics{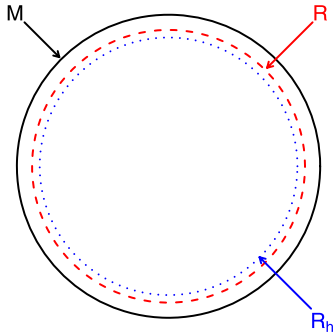}

\caption{The outer circle denotes the manifold $M$.
The dashed circle is the ridge $R$ of the density~$p$. The ridge is a
biased version of $M$
and acts as a surrogate for $M$.
The inner circle $R_h$ shows the ridge from a density estimator
with bandwidth $h$.
$R$ can be estimated at a much faster rate than $M$.}
\label{fig::surrogate-crop}
\end{figure}

\begin{example*}
A stylized example is shown in Figure~\ref{fig::Pedagogical}.
In this example, the density is
$p(x) = \int_M \phi(x-z) w(z) \,dz$
where $x\in\mathbb{R}^2$,
$M$ is a circle in $\mathbb{R}^2$,
$w$ is a smooth (but nonuniform) density supported on $M$
and $\phi$ is a two-dimensional Gaussian with a variance $\sigma^2$
that is much smaller than the radius of the circle.
The ridge $R$ is a
one-dimensional subset of $\mathbb{R}^2$.
The figure has a solid curve to show the ridge lifted onto $p$, that is,
the curve shows the set
$\{(x,p(x))\dvtx x\in R\}$.
The ridge $R$
does not coincide exactly with $M$
due to the blurring by convolution with the Gaussian.
In fact, $R$ is a circle with
slightly smaller radius than $M$.
That is, $R$ is a biased version of $M$.
Figure~\ref{fig::surrogate-crop} shows $M$ and $R$.

Note that the density is not uniform over the ridge.
Indeed, there can be modes
($0$-dimensional ridges) within a ridge.
What matters is that the function
rises sharply as we approach the ridge
(strongly negative eigenvalue). 
\end{example*}

One of the main points of this paper is that
$R$ captures the essential features of~$M$.
If we can live with the slight bias in $R$,
then it is better to estimate $R$ since $R$ can be estimated
at a polynomial rate while $M$ can only be estimated at a logarithmic rate.
Throughout this paper, we take the dimension
of interest $d$ as fixed and given.

Many different and useful definitions of a ``ridge'' have been proposed;
see the discussion of related work at the end of this section.
We make no claim as to the uniqueness and optimality of ours.
Our definition is motivated by four useful properties that we
demonstrate in this paper:
\begin{longlist}[1.]
\item[1.] If $\hat p$ is close to $p$, then $\hat R$ is close to $R$
where $\hat R$ is the ridge of $\hat p$ and $R$ is the ridge of $p$.
\item[2.] If the data-generating distribution is concentrated near a
manifold $M$, then the
ridge $R$ approximates $M$ both geometrically and topologically.
\item[3.]$R$ can be estimated at a polynomial rate, even in cases where
$M$ can
be estimated at only a logarithmic rate.
\item[4.] The definition corresponds essentially with the algorithm derived by
\citet{Principal}. That is, our definition provides
a mathematical formalization of their algorithm.
\end{longlist}

Our broad goal is to provide a theoretical framework for understanding
the problem of estimating hyper-ridge sets.
In particular, we show that the ridges of a kernel
density estimator consistently estimate the ridges of the density,
and we find and upper bound on the rate of convergence.
The main results of this paper are (stated here informally):
\begin{itemize}
\item\emph{Stability} (Theorem~\ref{theorem::close}).
If two densities are sufficiently close together, their hyper-ridge
sets are also close together.
\item\emph{Estimation} (Theorem~\ref{thm::main1}).
There is an estimator $\hat R$ such that
%
\begin{equation}
\operatorname{\mathsf{Haus}}(R,\hat R) = O_P \biggl( \biggl(\frac{\log n}{n}
\biggr)^{{2}/{(D+8)}} \biggr),
\end{equation}
where
$\operatorname{\mathsf{Haus}}$ is the Hausdorff distance, defined in equation (\ref
{eq::hausdorff}).
Moreover, $\hat R$~is topologically similar to $R$
in the sense that small dilations of these sets are topologically similar.
\item\emph{Surrogate} (Theorem~\ref{theorem::main}).
In the Hidden Manifold case with small noise variance $\sigma^2$
and assuming $M$ has no boundary,
the hyper-ridge set of the density $p$ satisfies
%
\begin{equation}
\operatorname{\mathsf{Haus}}(M,R) = O \bigl( \sigma^2 \log(1/\sigma) \bigr)
\end{equation}
and $R$ is topologically similar to $M$.
Hence, when the noise $\sigma$ is small,
the ridge is close to $M$.
Note that we treat $M$ as fixed while $\sigma\to0$.
It then follows that
%
\begin{equation}
\operatorname{\mathsf{Haus}}(M,\hat R) = O_P \biggl( \biggl(\frac{\log n}{n}
\biggr)^{{2}/{(D+8)}} \biggr) + O \bigl( \sigma^2 \log(1/\sigma)
\bigr).
\end{equation}
\end{itemize}

This leaves open the question of how to
locate the ridges of the density estimator.
Fortunately, this latter problem has recently been solved by
\citet{Principal}
who derived a practical algorithm
called the \emph{subspace constrained mean shift (SCMS) algorithm}
for locating the ridges.
\citet{Principal} derived their method assuming that
the underlying density function is known
(i.e., they did not discuss the effect of estimation error).
We, instead, assume the density is estimated from a finite sample
and adapt their algorithm accordingly by including a denoising step
in which we discard points with low density.
This paper provides a
statistical justification for, and extension to, their algorithm.
We introduce a modification of their algorithm called
SuRF (Subspace Ridge Finder) that applies
density estimation, followed by denoising, followed by SCMS.

\emph{Related work.}
Zero dimensional ridges are modes
and in this case
ridge finding reduces to mode estimation
and SCMS reduces to the
mean shift clustering algorithm
[\citet{Fukunaga,cheng1995mean,li2007nonparametric,chacon}].

If the hidden structure is a manifold,
then the process of finding the structure is known
as \emph{manifold estimation} or
\emph{manifold learning}.
There is a large literature on manifold estimation and related techniques.
Some useful references are
\citet{smale}
\citet{CCDM},
\citeauthor{us::2010}
(\citeyear{genovese2009path,genovese2012geometry,singular::2011,us::2010}),
\citet{tenenbaum2000global,roweis2000nonlinear}
and references therein.

The notion of ridge finding spans many fields.
Previous work on ridge finding
in the statistics literature includes
\citet{cheng2004estimating},
\citet{hall2001local},
\citet{wegman2002smoothings,wegman2} and
\citet{hall1992ridge}.
These papers focus on visualization and exploratory analysis.
An issue that has been discussed extensively
in the applied math and computer science literature is
how to define a ridge.
A detailed history and taxonomy is given in the text by
\citet{eberly1996ridges}.
Two important classes of ridges are watershed ridges, which are global
in nature,
and height ridges, which are locally defined.
There is some debate about
the virtues of various definitions.
See, for example,
\citet{norgard2012second,norgard2012second2}.
Related definitions also appear in the fluid dynamics literature
[\citet{schindler2012ridge}] and
astronomy
[\citet{aragon2010spine,sousbie2008three}].
There is also a literature on Reeb graphs
[\citet{ge2011data}]
and metric graphs
[\citet{aanjaneya2012metric,lecci2013statistical}].
Metric graph methods are ideal for representing intersecting filamentary
structure but are much more sensitive to noise than the methods in this paper.
It is not our intent in this paper to argue
that one particular definition of ridge is optimal for all purposes.
Rather, we use a particular definition which is well suited for
studying the statistical estimation of ridges.

More generally, there is a vast literature on
hunting for structure in point clouds and analyzing the shapes
of densities.
Without attempting to be exhaustive,
some representative work includes
\citet{davenport2010joint},
\citet{klemela2009smoothing,adams2011morse,chazal2011persistence,bendich2012local}.

Throughout the paper,
we use symbols like
$C, C_0,C_1, c,c_0, c_1,\ldots$ to denote generic positive constants
whose value may be different in different expressions.

\section{Model and ridges}
\label{sec::conditions}

In this section, we describe our assumptions about the
data and give a formal definition of hyper-ridge sets,
which we call \emph{ridges} from now on.
Further properties of ridges are stated and proved
in Section~\ref{section::ridges}.

We start with a point cloud $X_1,\ldots,X_n\in\R^D$.
We assume that these data
comprise a random sample from a distribution $P$ with density $p$,
where $p$ has at least five bounded, continuous derivatives.
This is all we assume for the density ridge case.
In the hidden manifold case, we assume further that $P$ and $p$ are
derived from a $d$-dimensional manifold $M$
by convolution with a noise distribution,
where $d < D$.
Specifically,
we assume that $M$ is embedded within a compact subset $\cK\subset\R^D$
and that
%
\begin{equation}
\label{eq::model}
P = (1-\eta) \operatorname{Unif}(\cK) + \eta(W \star\Phi_\sigma),
\end{equation}
where
$0 < \eta\le1$,
$\operatorname{Unif}(\cK)$ is a uniform distribution on $\cK$,
$\star$ denotes convolution,
$W$ is a distribution supported on $M$,
and
$\Phi_\sigma$ is a Gaussian distribution on $\R^D$ with
zero mean and covariance $\sigma I_D$.
While we could consider a more general noise distribution in (\ref{eq::model}),
we focus on the common assumption of Gaussian noise.
In that case,
a hidden manifold $M$ can only be estimated at a logarithmic rate
[\citet{singular::2011}],
so ridge estimators are particularly valuable.
(Even when $M$ can be estimated at a polynomial rate, ridge estimators are
often easier in practice than estimating the manifold,
which would involve deconvolution.)

The data generating process under model (\ref{eq::model})
is equivalent to the following steps:
\begin{longlist}[1.]
\item[1.] Draw $B$ from a $\operatorname{Bernoulli}(\eta)$.
\item[2.] If $B=0$, draw $X$ from a uniform distribution on $\cK$.
\item[3.] If $B=1$, let $X= Z + \sigma\varepsilon$ where
$Z\sim W$ and $\varepsilon$ is additional noise.
\end{longlist}
Points $X_i$ drawn from $\operatorname{Unif}(\cK)$ represent background clutter.
Points $X_i$ drawn from
$W\star\Phi_\sigma$
are noisy observations from $M$.
When $M$ consists of a finite set of points,
this can be thought of as a clustering model.

\subsection{Definition of ridges}
As in \citet{Principal}, our definition of
ridges relies on the gradient and Hessian
of the density function $p$.
Recall that $0 < d < D$ is fixed throughout.
Given a function $p\dvtx \R^D \to\R$, let
$g(x) = \nabla p(x)$
denote its gradient and
$H(x)$ its Hessian matrix,
at $x$.
Let
%
\begin{equation}
\lambda_1(x) \ge\lambda_2(x) \ge\cdots\ge
\lambda_{d}(x) > \lambda _{d+1}(x) \ge\cdots\ge
\lambda_D(x)
\end{equation}
denote the eigenvalues of $H(x)$
and let $\Lambda(x)$ be the diagonal matrix whose
diagonal elements are the eigenvalues.
Write the spectral decomposition of $H(x)$ as
$H(x) = U(x) \Lambda(x) U(x)^T$.
Let $V(x)$ be the last $D-d$ columns of $U(x)$
(i.e., the columns corresponding to the $D-d$ smallest eigenvalues).
If we write $U(x) = [V_\diamond(x)\dvtx V(x)]$ then we can write
$H(x) = [V_\diamond(x)\dvtx V(x)]\Lambda(x) [V_\diamond(x)\dvtx\break  V(x)]^T$.
Let $L(x) \equiv L(H(x))= V(x) V(x)^T$ be the projector
onto the linear space defined by
the columns of $V(x)$.
We call this the local normal space
and the space spanned by
$L^\perp(x) = I- L(x) = V_\diamond(x) V_\diamond(x)^T$
is the local tangent space.
Define the \emph{projected gradient}
%
\begin{equation}
G(x) = L(x) g(x).
\end{equation}

If the vector field $G(x)$ is Lipschitz
then by Theorem~3.39 of \citet{irwin1980smooth},
$G$~defines a global flow as follows.
The flow is a family of functions
$\phi(x,t)$ such that
$\phi(x,0)=x$ and $\phi'(x,0) = G(x)$ and
$\phi(x,s+t) = \phi(\phi(x,t),s)$.
The flow lines, or integral curves, partition the space
(see Lemma~\ref{lemma::unique})
and at each $x$ where $G(x)$ is nonnull,
there is a unique integral curve passing through $x$.
Thus, there is one and only one flow line through each nonridge point.
The intuition is that the flow passing through $x$ is
a gradient ascent path moving toward higher values of~$p$.
Unlike the paths defined by the gradient $g$ which move toward modes,
the paths defined by
the projected gradient $G$ move toward ridges.
The SCMS algorithm, which we describe later,
can be thought of
as approximating the flow with
discrete, linear steps $x_{k+1}\leftarrow x_k + h G(x_k)$.
[A proof that the linear interpolation of these points
approximates the flow in the case $d=0$
is given in \citet{arias2013estimation}.]

A map $\pi\dvtx  \R\to\R^D$
is an integral curve with respect to the flow of $G$ if
%
\begin{equation}
\label{eq::diffeq} \pi'(t) = G\bigl(\pi(t)\bigr) = L\bigl(\pi(t)\bigr) g
\bigl(\pi(t)\bigr).
\end{equation}

\emph{Definition}:
The \emph{ridge $R$} of dimension $d$ is given by
$R = \{x\dvtx \Vert G(x)\Vert =0, \lambda_{d+1}(x) < 0\}$.\vadjust{\goodbreak}

Note that the ridge consists of the destinations of the integral curves:
$y\in R$ if
$\lim_{t\to\infty} \pi(t) = y$ for some $\pi$
satisfying (\ref{eq::diffeq}).

Our definition is motivated by
\citet{Principal}
but is slightly different.
They first define the $d$-critical points as those for which
$\Vert G(x)\Vert  =0$.
They call a critical point regular if it is $d$-critical but not
$(d-1)$-critical.
Thus, a~mode within a one-dimensional ridge
is not regular.
A regular point with $\lambda_{d+1} < 0$ is called
a principal point.
According to our definition, the ridge lies between the critical set
and the principal set.
Thus, if a mode lies on a one-dimensional ridge,
we include that point as part of the ridge.


\subsection{Assumptions}

We now record the main assumptions about the ridges that we will
require for the results.

\emph{Assumption \textup{(A0)} differentiability.}
For all $x$,
$g(x)$, $H(x)$ and $H'(x)$ exist.

\emph{Assumption \textup{(A1)} eigengap.}
Let $B_D(x,\delta)$ denote a $D$-dimensional ball of radius $\delta$
centered at $x$ and let
$R\oplus\delta= \bigcup_{x\in R} B_D(x,\delta)$.
We assume that there exists
$\beta>0$ and $\delta>0$
such that, for all $x\in R\oplus\delta$,
$\lambda_{d+1}(x) < -\beta$ and
$\lambda_d(x) - \lambda_{d+1}(x) > \beta$.

\emph{Assumption \textup{(A2)} path smoothness.}
For each $x\in R\oplus\delta$,
{\renewcommand{\theequation}{A2}
\begin{equation}\label{eqa2}
\bigl\Vert L^\perp(x)g(x)\bigr\Vert \bigl\Vert H'(x)\bigr\Vert
_{\mathrm{max}} < \frac{\beta^2}{2 D^{3/2}},
\end{equation}}
\hspace*{-2pt}where
$H'(x) = \frac{d \operatorname{\mathsf{vec}}(H(x))}{d x^T}$,
$L^\perp= I-L$ and
$\Vert A\Vert _{\mathrm{max}} =\max_{j,k} |A_{jk}|$.

Condition (A1) says that $p$ is
sharply curved around the ridge
in the $D-d$ dimensional space normal to the ridge.
To give more intuition about the condition,
consider the problem of estimating a mode in one dimension.
At a mode $x$, we have that $p'(x) =0$ and
$p''(x) < 0$.
However, the mode cannot be uniformly consistently estimated
by only requiring the second derivative to be negative
since $p''(x)$ could be arbitrarily close to 0.
Instead, one needs to assume that
$p''(x) < -\beta$ for some positive constant $\beta$.
Condition (A1) may be thought of as
the analogous condition for a ridge.
(\ref{eqa2}) is a third derivative condition which
implies that the paths cannot be too wiggly.
(\ref{eqa2}) also constrains the gradient from being too steep in the
perpendicular direction.
Note that these conditions are local: they hold in a size $\delta$ neighborhood
around the ridge.

\section{Technical background}
Now we review some background.
We recommend that the reader quickly skim this section and then refer
back to it
as needed.

\subsection{Distance function and Hausdorff distance}

We let $B(x,r)\equiv B_D(x,r)$
denote a $D$-dimensional open ball centered at $x\in\mathbb{R}^D$ with
radius $r$.
If $A$ is a set and $x$ is a point then we define the
\emph{distance function}
%
\setcounter{equation}{7}
\begin{equation}
\label{eq::distfun} d_A(x) =d(x,A) = \inf_{y\in A}\Vert x-y
\Vert,\vadjust{\goodbreak}
\end{equation}
where
$\Vert  \cdot\Vert $ is the Euclidean norm.
Given two sets $A$ and $B$, the \emph{Hausdorff distance} between $A$
and $B$ is
%
\begin{equation}
\operatorname{\mathsf{Haus}}(A,B) = \inf \{ \varepsilon\dvtx A\subset B\oplus\varepsilon
\mbox{ and }
B\subset A\oplus\varepsilon \}=
\sup_x \bigl|d_A(x) - d_B(x)\bigr|, \label{eq::hausdorff}
\end{equation}
where
%
\begin{equation}
A\oplus\varepsilon= \bigcup_{x\in A} B_D(x,
\varepsilon) = \bigl\{x\dvtx d_A(x) \leq\varepsilon \bigr\}
\end{equation}
is called the \emph{$\varepsilon$-dilation} of $A$.
The dilation can be thought of as a smoothed version of $A$.
For example, if there are any small holes in $A$, these will be filled in
by forming the dilation $A\oplus\varepsilon$.

We use Hausdorff distance to measure the distance between sets
for several reasons: it is the most commonly used distance between sets,
it is a very strict distance and
is analogous to the familiar $L_\infty$ distance
between functions for sets.

\subsection{Topological concepts}

This subsection follows
Chazal, Cohen-Steiner and Lieutier (\citeyear{chazal2009sampling})
and \citet{chazal2005lambda}.
The \emph{reach} of a set $K$,
denoted by $\operatorname{\mathsf{reach}}(K)$,
is the largest $r>0$ such that each point in
$K\oplus r$ has a unique projection onto $K$.
A set with positive reach is, in a sense, a smooth set
without self-intersections.

Now we describe a generalization of reach
called $\mu$-reach.
The key point is simply that the $\mu$-reach is weaker than reach.
The full details can be found in the aforementioned references.
Let $A$ be a compact set.
Following
Chazal and Lieutier (\citeyear{chazal2005lambda}) define the gradient
$\nabla_A(x)$ of $d_A(x)$ to be the usual
gradient function
whenever this is well defined.
However, there may be points $x$ at which $d_A$ is not
differentiable in the usual sense.
In that case,
define the gradient as follows.
For $x\in A$ define $\nabla_A(x)=0$ for all $x\in A$.
For $x\notin A$, let
$\Gamma(x) = \{y\in A\dvtx \Vert x-y\Vert  = d_A(x)\}$.
Let $\Theta(x)$ be the center of the unique smallest closed ball
containing $\Gamma(x)$.
Define
$\nabla_A(x) = \frac{x - \Theta(x)}{d_A(x)}$.\vspace*{1pt}

The \emph{critical points}
are the points at which $\nabla_A(x)=0$.
The \emph{weak feature size}
$\operatorname{\mathsf{wfs}}(A)$ is the distance from $A$ to its
closest critical point.
For $0 < \mu< 1$, the \emph{$\mu$-reach} $\operatorname{\mathsf{reach}}_\mu(A)$ is
$\operatorname{\mathsf{reach}}_\mu(A) = \inf\{d\dvtx \chi(d) < \mu\}$
where
$\chi(d) = \inf\{ \Vert \nabla_A (x)\Vert \dvtx d_A(x)=d\}$.
It can be shown that
$\operatorname{\mathsf{reach}}_\mu$ is nonincreasing in $\mu$,
that
$\operatorname{\mathsf{wfs}}(A) = \lim_{\mu\to0}\operatorname{\mathsf{reach}}_\mu(A)$ and that
$\operatorname{\mathsf{reach}}(A)=\lim_{\mu\to1}\operatorname{\mathsf{reach}}_\mu(A)$.

\begin{figure}

\includegraphics{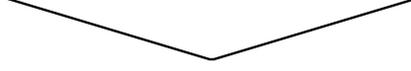}

\caption{A straight line as infinite reach.
A line with a corner, as in this figure,
has 0 reach but has positive $\mu$-reach.}
\label{fig::mureach}
\end{figure}

As a simple example,
a circle $C$ with radius $r$ has
$\operatorname{\mathsf{reach}}(R)=r$.
However, if we bend the circle slightly to create a corner,
the reach is 0 but, provided the kink is not too extreme,
the $\mu$-reach is still positive.
As another example, a straight line as infinite reach.
Now suppose we add a corner as in Figure~\ref{fig::mureach}.
This set has 0 reach but has positive $\mu$-reach.

Two maps
$f\dvtx A \to B$ and
$g\dvtx A \to B$ are \emph{homotopic} if there exists a continuous map
$H\dvtx [0,1]\times A \to B$ such that
$H(0,x)=f(x)$ and
$H(1,x)=g(x)$.
Two sets $A$ and $B$ are
homotopy equivalent if there are continuous maps
$f\dvtx A \to B$ and $g\dvtx B \to A$
such that the following is true:
(i) $g\circ f$ is homotopic to the identity map on $A$ and
(ii) $f\circ g$ is homotopic to the identity map on $B$.
In this case we write
$A\cong B$.
Sometimes $A$ fails to be homotopic to $B$
but $A$ is homotopic to
$B\oplus\delta$ for every sufficiently small $\delta>0$.
This happens because
$B\oplus\delta$ is slightly smoother than $B$.
If $A \cong B\oplus\delta$ for all small $\delta>0$,
we will say that $A$ and $B$ are \emph{nearly homotopic} and we will write
$A\VHF B$.

The following result
[Theorem~4.6 in Chazal, Cohen-Steiner and Lieutier
(\citeyear{chazal2009sampling})]
says that if a set $K$ is smooth and
$\tilde{K}$ is close to $K$,
then a smoothed version of $\tilde{K}$ is nearly homotopy equivalent
to $K$.

\begin{theorem}[{[Chazal, Cohen-Steiner and Lieutier (\citeyear{chazal2009sampling})]}]
\label{thm::CCSL}
Let $K$ and $\tilde{K}$ be compact sets and let
$\varepsilon= \operatorname{\mathsf{Haus}}(\tilde{K},K)$.
If
%
\begin{equation}
\varepsilon< \frac{\mu^2 \operatorname{\mathsf{reach}}_\mu(K)}{5\mu^2 + 12}
\quad\mbox{and}\quad \frac{4 \varepsilon}{\mu^2} \leq\alpha< \operatorname{\mathsf{reach}}_\mu(K) - 3\varepsilon
\end{equation}
then
$(\tilde{K}\oplus\alpha)  \VHF K$.
\end{theorem}

\subsection{Matrix theory}

We make extensive use of matrix theory
as can be found in
\citet{stewart1990matrix},
\citet{bhatia1997matrix},
\citet{Horn} and
\citet{magnus1988matrix}.

Let $A$ be an $m\times n$ matrix.
Let $A_{jk}$ denote an element of the matrix.
Then the Frobenius norm is
$\Vert A\Vert _F = \sqrt{\sum_{j,k} A_{jk}^2}$
and the operator norm is
$\Vert A\Vert  = \sup_{\Vert x\Vert =1} \Vert Ax\Vert $.
We define
$\Vert A\Vert _{\mathrm{max}} = \max_{j,k}|A_{jk}|$.
It is well known that
$\Vert A\Vert \leq\Vert A\Vert _F \leq\sqrt{n}\Vert A\Vert $,
that
$\Vert A\Vert _{\mathrm{max}} \leq\Vert A\Vert  \leq\sqrt{mn} \Vert A\Vert _{\mathrm{max}}$
and that
$\Vert A\Vert _F \leq\sqrt{mn}\Vert A\Vert _{\mathrm{max}}$.

The $\operatorname{\mathsf{vec}}$ operator
converts a matrix into a vector
by stacking the columns.
Thus, if $A$ is $m\times n$ then
$\operatorname{\mathsf{vec}}(A)$ is a vector of length $mn$.
Conversely, given a vector $a$ of length $mn$,
let $[[a]]$ denote the $m\times n$ matrix obtained by stacking $a$
columnwise into matrix form.
We can think of $[[a]]$ as the ``anti-vec'' operator.

If $A$ is $m\times n$ and $B$ is $p\times q$ then
the Kronecker $A\otimes B$ is the $mp\times nq$ matrix
%
\begin{equation}
\left[ \matrix{ A_{11} B & \cdots&
A_{1n}B
\vspace*{2pt}\cr
\vdots& & \vdots
\vspace*{2pt}\cr
A_{m1}B & \cdots& A_{mn}B }
 \right].
\end{equation}
If $A$ and $B$ have the same dimensions,
then the Hadamard product
$C=A\circ B$ is defined by $C_{jk} = A_{jk} B_{jk}$.

For matrix calculus, we follow the conventions in
\citet{magnus1988matrix}.
If $F\dvtx \R^D\to\R^k$ is a vector-valued map
then the Jacobian matrix
will be denoted by
$F'(x)$ or $dF/dx$.
This is the $D\times k$ matrix with
$F'(x)_{jk} = \partial F_i(x)/\partial x_j$.
If $F\dvtx \R^D\to\R^{m\times p}$ is a matrix-valued map
then $F'(x)$ is a $mp\times D$ matrix defined by
%
\begin{equation}
\label{eq::Fprime} F'(x) \equiv\frac{dF}{dx^T}= \frac{d \operatorname{vec}(F(x))}{dx^T}.
\end{equation}
If
$F\dvtx \R^{n\times q} \to\R^{m\times p}$
then the derivative is a
$mp\times nq$ matrix given by
\[
F'(X) \equiv\frac{dF}{dX} = \frac{d\operatorname{vec}(F(X))}{d\operatorname{vec}(X)^T}.
\]
We then have the following product rule for matrix calculus:
if $F\dvtx \R^D \to\R^{m\times p}$ and
$G\dvtx \R^D \to\R^{p\times q}$
then
\[
\frac{d F(x)G(x)}{dx} = \bigl(G^T(x) \otimes I_m\bigr)
F'(x) + \bigl(I_q \otimes F(x)\bigr)
G'(x).
\]
Also, if $A(x) = f(x) I$ then
$A'(x) = \operatorname{\mathsf{vec}}(I) \otimes(\nabla f(x))^T$
where $\nabla f$ denotes the gradient of $f$.

The following version of the Davis--Kahan theorem is from
\citet{von2007tutorial}.
Let $H$ and $\tilde{H}$ be two symmetric, square $D\times D$ matrices.
Let $\Lambda$ be the diagonal matrix of eigenvalues of $H$.
Let $S\subset\R$ and let
$V$ be the matrix whose columns
are the eigenvectors corresponding\vspace*{1pt} to the
eigenvalues of $H$ in $S$
and similarly for
$\tilde V$ and $\tilde H$.
Let
%
\begin{equation}
\label{eq::beta} \beta= \min \bigl\{ |\lambda-s|\dvtx \lambda\in\Lambda\cap
S^c, s\in S \bigr\}.
\end{equation}
%
According
to the Davis--Kahan theorem,
%
\begin{equation}
\label{eq::davis-kahan} \bigl\Vert VV^T - \tilde V \tilde V^T\bigr\Vert \leq
\frac{\Vert H - \tilde H\Vert _F}{\beta}.
\end{equation}

Let $H$ be a $D\times D$ square, symmetric matrix with eigenvalues
$\lambda_1 \ge\cdots\ge\lambda_D$.
Let $\tilde H$ be another square, symmetric matrix with eigenvalues
$\tilde\lambda_1 \ge\cdots\ge\tilde\lambda_D$.
By Weyl's theorem
[Theorem~4.3.1 of \citet{Horn}],
we have that
%
\begin{equation}
\label{eq::weyl1} \lambda_n(\tilde H - H) + \lambda_i(H)
\leq\lambda_i(\tilde H) \leq \lambda_i(H) +
\lambda_1(\tilde H - H).
\end{equation}
It follows easily that
%
\begin{equation}
\label{eq::weyl2} \bigl|\lambda_i(H) - \lambda_i(\tilde H)\bigr|
\leq\Vert H-\tilde H\Vert \leq D \Vert H-\tilde H\Vert _{\mathrm{max}}.
\end{equation}

\section{Properties of ridges}
\label{section::ridges}

In this section, we examine some of the properties
of ridges as they were defined in Section~\ref{sec::conditions}
and show that, under appropriate conditions,
if two functions are close together then their ridges are close and
are topologically similar.

\subsection{Arclength parameterization}
It will be convenient to
parameterize the gradient ascent paths by arclength.
Thus, let $s\equiv s(t)$ be the arclength from
$\pi(t)$ to $\pi(\infty)$:
%
\begin{equation}
\label{eq::arclength} s(t)=\int_{t}^\infty\bigl\Vert
\pi'(u)\bigr\Vert \,du.
\end{equation}
Let $t\equiv t(s)$ denote the inverse of $s(t)$.
Note that
%
\begin{equation}\quad
t'(s) = - \frac{1}{\Vert \pi'(t(s))\Vert } = - \frac{1}{\Vert L(\pi(t(s)))g(\pi(t(s)))\Vert } = -
\frac{1}{\Vert G(\pi(t(s)))\Vert }.
\end{equation}
Let $\gamma(s) = \pi(t(s))$.
Then
%
\begin{equation}
\gamma'(s) = - \frac{G(\gamma(s))}{\Vert G(\gamma(s))\Vert },
\end{equation}
which is a restatement of (\ref{eq::diffeq}) in
the arclength parameterization.

In what follows,
we will often abbreviate notation
by using the subscript $s$ in the following way:
$G_s = G(\gamma(s)), H_s = H(\gamma(s)),\ldots,$
and so forth.

\subsection{Differentials}
\label{subsection::differentials}

We will need derivatives of $g$, $H$, and $L$.
The derivative of $g$ is the Hessian $H$.
Recall from (\ref{eq::Fprime}) that
$H'(x) = \frac{d \operatorname{\mathsf{vec}}(H(x))}{d x^T}$.
We also need derivatives along the curve $\gamma$.
The derivative of a functions $f$ along $\gamma$ is
%
\begin{equation}
\dt{f}_{\gamma(s)} \equiv\dt{f}_s = \lim_{\varepsilon\to0}
\frac{ f(\gamma(s+\varepsilon))- f(\gamma(s))}{\varepsilon}.
\end{equation}
Thus, the derivative of the gradient $g$ along $\gamma$ is
%
\begin{equation}
\dt{g}_{\gamma(s)} \equiv\dt{g}_s = \lim_{\varepsilon\to0}
\frac{ g(\gamma(s+\varepsilon))- g(\gamma(s))}{\varepsilon} = H_s \gamma _s' = -
\frac{H_s G_s}{\Vert G_s\Vert }.
\end{equation}

We will also need the derivative of $H$
in the direction of a vector $z$ which we will denote by
\[
H'(x;z) \equiv \lim_{\varepsilon\to0} \frac{ H(x+\varepsilon z) - H(x)}{\varepsilon}.
\]
We can write an explicitly formula for
$H'(x;z)$ as follows.
Note that the elements of $H'$ are the
partial derivatives
$\partial H_{jk}(x)/\partial x_\ell$
arranged in a $D^2 \times D$ matrix.
Hence, $H'(x;z) = [[H'(x) z]]$.
(Recall that $[[a]]$ stacks a vector into a matrix.)
Note that
$[[H'(x) z]]$ is a $D\times D$ matrix.

Recall that $L(x) \equiv L(H(x)) = V(x) V(x)^T$.
The collection
$\{L(x)\dvtx x\in\R^D\}$
defines a matrix field: there is a matrix $L(x)$
attached to each point $x$.
We will need the derivative of this field along the integral curves
$\gamma$.
For any $x\notin R$, there is a unique path $\gamma$ and unique $s>0$
such that
$x=\gamma(s)$.
Define
%
\begin{eqnarray}
\dt{L}_s &\equiv&\dt{L}(x) \equiv \lim_{\varepsilon\to0}
\frac{ L(H(\gamma(s+\varepsilon))) - L(H(\gamma
(s)))}{\varepsilon}
\nonumber
\\[-8pt]
\\[-8pt]
\nonumber
&=&\lim_{t\to0}\frac{ L(H + t E) - L(H)}{t},
\end{eqnarray}
where $H = H(\gamma(s))$
and $E= (d/ds)H(\gamma(s)) = H'(x;z)$ with $z=\gamma'(s)$.

\subsection{Uniqueness of the \texorpdfstring{$\gamma$}{gamma} paths}

\begin{lemma}\label{lemma::unique}
Conditions \textup{(A0)--(\ref{eqa2})} imply that,
for each $x\in(R\oplus\delta) - R$,
there is a unique path $\gamma$ passing through $x$.
\end{lemma}

\begin{pf}
We will show that
the vector field $G(x)$ is Lipschitz
over $R\oplus\delta$.
The result then follows from
Theorem~3.39 of \citet{irwin1980smooth}.
Recall that $G=Lg$ and $g$ is differentiable.
It suffices to show that
$L$ is differentiable over $R\oplus\delta$.
Now $L(x) = L(H(x))$.
It may be shown that,
as a function of $H$,
$L$ is Frechet differentiable.
And $H$ is differentiable by assumption.
By the chain rule, $L$ is differentiable
as a function of $x$.
Indeed,
$dL/dx$ is the $D^2\times D$ matrix
whose $j$th column is
$\operatorname{\mathsf{vec}}(L^\dagger E_j)$
where $E_j = [[H' e_j]]$,
$L^\dagger$ denotes the Frechet derivative,
and $e_j$ is the vector which is 1 in the
$j$th coordinate and zero otherwise.
\end{pf}

\subsection{Quadratic behavior}

Conditions (A1) and (\ref{eqa2}) imply that the function $p$
has quadratic-like behavior near the ridges.
This property is needed for establishing the convergence of
ridge estimators.
In this section, we formalize this notion of quadratic behavior.
Give a path $\gamma$,
define the function
%
\begin{equation}
\xi(s) = p\bigl(\pi(\infty)\bigr) - p\bigl(\pi\bigl(t(s)\bigr)\bigr) = p\bigl(
\gamma(0)\bigr) - p\bigl(\gamma(s)\bigr).
\end{equation}
Thus, $\xi$ is simply the drop in the function $p$ along the curve
$\gamma$
as we move away from the ridge.
We write $\xi_x(s)$ if we want to emphasize that
$\xi$ corresponds to the path $\gamma_x$ passing through the point $x$.
Since $\xi\dvtx [0,\infty) \to[0,\infty)$,
we define its derivatives in the usual way, that is,
$\xi'(s) = d\xi(s)/ds$.

\begin{lemma}
\label{lemma::thiss}
Suppose that \textup{(A0)--(\ref{eqa2})} hold.
For all $x\in R\oplus\delta$, the following are true:
\begin{longlist}[1.]
\item[1.]$\xi(0) = 0$.
\item[2.]$\xi'(s) = \Vert G(\gamma(s))\Vert $ and $\xi'(0)=0$.
\item[3.] The second derivative of $\xi$ is:
%
\begin{equation}
\xi''(s) = -\frac{G_s^T H_s G_s}{\Vert G_s\Vert ^2} + \frac{g_s^T \dt{L}_s G_s}{\Vert G_s\Vert }.
\end{equation}
\item[4.]$\xi''(s) \ge\beta/2$.
\item[5.]$\xi(s)$ is nonincreasing in $s$.
\item[6.]$\xi(s) \ge\frac{\beta}{4}\Vert \gamma(0)-\gamma(s)\Vert ^2$.
\end{longlist}
\end{lemma}

\begin{pf}
1. The first condition $\xi(0) =0$ is immediate from the definition.

2. Next,
\begin{eqnarray*}
\xi'(s) &=& - \frac{d p(\gamma(s))}{ds} = -g_s
\gamma_s' = \frac{g_s^T G_s}{\Vert G_s\Vert } = \frac{g_s^T L_s g_s}{\Vert G_s\Vert }
\\
&=& \frac{g_s^T L_s L_s g_s}{\Vert G_s\Vert } =
\frac{G_s^T G_s}{\Vert G_s\Vert } = \Vert G_s\Vert.
\end{eqnarray*}
Since the projected gradient is 0 at the ridge,
we have that
$\xi'(0)=0$.

3.
Note that
$(\xi'(s))^2 = \Vert G_s\Vert ^2 =
G_s^T G_s = g_s^T L_s g_s \equiv a(s)$.
Differentiating both sides of this equation, we have that
$2\xi'(s) \xi''(s) = a'(s)$, and hence
\[
\xi''(s) = \frac{a'(s) }{2 \xi'(s)} = \frac{a'(s)}{2 \Vert G_s\Vert }.
\]
Now
%
\begin{equation}
\label{eq::as} a'(s) = (\dt{g}_s)^T
L_s g_s + g_s^T
\dt{L}_s g_s + g_s^T
L_s \dt {g}_s = 2(\dt{g}_s)^T
L_s g_s + g_s^T
\dt{L}_s g_s.
\end{equation}
Since
$L_s L_s = L_s$ we have that
$\dt{L}_s = L_s \dt{L_s} + \dt{L}_s L_s$,
and hence
\begin{eqnarray*}
g_s^T \dt{L}_s g_s &=
g_s^T L_s \dt{L}_s
g_s + g_s^T \dt{L}_s
L_s g_s = G_s^T
\dt{L}_s g_s + g_s^T
\dt{L}_s G_s = 2 g_s^T
\dt{L}_s G_s.
\end{eqnarray*}
Therefore,
%
\begin{equation}
\label{eq::as2} a'(s)= 2(\dt{g}_s)^T
L_s g_s + 2 g_s^T
\dt{L}_s G_s.
\end{equation}
Recall that
$\dt{g}_s = -\frac{H_s G_s}{\Vert G_s\Vert }$.
Thus,
%
\begin{equation}
\xi''(s) =\frac{a'(s)}{2\Vert G_s\Vert } = -\frac{G_s^T H_s G_s}{\Vert G_s\Vert ^2} +
\frac{g_s^T \dt{L}_s G_s}{\Vert G_s\Vert }.
\end{equation}

4. The first term in $\xi''(s)$ is
$-\frac{G_s^T H_s G_s}{\Vert G_s\Vert ^2}$.
Since $G$ is in the column space of $V$,
$G_s^T H_s G_s = G_s^T (V_s \Lambda_s V_s^T) G_s$
where
$\Lambda_s = \operatorname{diag}(\lambda_{d+1}(\gamma(s)),\ldots,\lambda
_D(\gamma(s)))$.
Hence, from (A1),
\[
\frac{G_s^T H_s G_s}{\Vert G_s\Vert ^2} = \frac{G_s^T (V_s \Lambda_s V_s^T) G_s}{\Vert G_s\Vert ^2} \leq \lambda_{\mathrm{max}}
\bigl(V_s \Lambda_s V_s^T\bigr)
< -\beta
\]
and thus
\[
-\frac{G_s^T H_s G_s}{\Vert G_s\Vert ^2} \ge\beta.
\]


Now we bound the second term
$\frac{g_s^T \dt{L}_s G_s}{\Vert G_s\Vert }$.
Since
$L_s + L^\perp_s = I$ and
$L_s G_s = G_s$, we have
$g_s^T \dt{L}_s G_s =
g_s^T L_s \dt{L}_s G_s + g_s^T L_s^\perp\dt{L}_s G_s=
g_s^T L_s \dt{L}_s L_s G_s + g_s^T L_s^\perp\dt{L}_s L_s G_s$.
Now
$|g_s^T L_s \dt{L}_s L_s G_s| = 0$.
To see this, note that
$L_s L_s = L_s$ implies
$L_s \dt{L}_s + \dt{L}_s L_s = \dt{L}_s$ implies
$L_s \dt{L}_s L_s + \dt{L}_s L_s = \dt{L}_sL_s$ implies
$L_s \dt{L}_s L_s =0$.
To bound $g_s^T L_s^\perp\dt{L}_s L_s G_s$ we proceed as follows.
Let $E = (d/ds) H(\pi(\gamma(s))) = H'(x;z)$ with $z=\gamma'(s)$.
Then, from Davis--Kahan,
\begin{eqnarray*}
\bigl|g_s^T L_s^\perp
\dt{L}_s L_s G_s\bigr| &=& \lim
_{t\to0} \frac{|g_s^T L_s^\perp(L(H+ tE) - L(H)) L_s G_s}{t}
\\
& \leq& \bigl\Vert L^\perp g_s\bigr\Vert \lim_{t\to0}
\frac{\Vert  (L(H+ tE) - L(H))\Vert }{t} \Vert G_s\Vert \\
&\leq& \frac{\Vert L^\perp g_s\Vert  \Vert E\Vert  \Vert G_s\Vert }{\beta}.
\end{eqnarray*}
Note that
$\Vert H'(x;z)\Vert  \leq D \Vert H'(x;z)\Vert _{\mathrm{max}} \leq D^{3/2}
 \Vert H'(x)\Vert _{\mathrm{max}} \Vert z\Vert  = D^{3/2} \times\break \Vert H'(x)\Vert _{\mathrm{max}}$.
So
$\frac{|g_s^T \dt{L}_s G_s|}{\Vert G_s\Vert } \leq D^{3/2} \Vert L^\perp g\Vert  \Vert H'\Vert _{\mathrm{max}}/\beta$
which is less than $\beta/2$ by~(\ref{eqa2}).
Therefore,
$\xi''(s) \geq\beta- (\beta/2) = \beta/2$.

5. Follows from 2.

6. For some $0\leq\tilde{s} \leq s$,
\[
\xi(s) = \xi(0) + s \xi'(0) + \frac{s^2}{2}
\xi''(\tilde{s}) = \frac{s^2}{2}
\xi''(\tilde{s}) \ge\frac{\beta s^2}{4}
\]
from part (4). So
\[
\xi(s)-\xi(0) \ge\frac{\beta}{4}s^2 \ge\frac{\beta}{4}\bigl\Vert
\gamma (0) - \gamma(s)\bigr\Vert ^2.
\]
\upqed\end{pf}

\subsection{Stability of ridges}

We now show that if two functions
$p$ and $\tilde{p}$ are close,
then their corresponding ridges
$R$ and $\tilde{R}$ are close.
We use $\tilde{g}, \tilde{H},\ldots$ etc.
to refer to the gradient, Hessian and so on, defined by $\tilde{p}$.
For any function $f\dvtx \R^D \to\R$, let
$\Vert f\Vert _\infty= \sup_{x\in R\oplus\delta}|f(x)|$.
Let
%
\begin{eqnarray}
\varepsilon&=& \Vert p-\tilde p\Vert _\infty,\qquad\varepsilon' =
\max_j \Vert g_j-\tilde g_j\Vert
_\infty,
\\
\varepsilon'' &=& \max_{jk}\Vert
H_{jk}-\tilde H_{jk}\Vert _\infty,\qquad
\varepsilon''' = \max
_{jk}\bigl\Vert H_{jk}'-\tilde
H_{jk}'\bigr\Vert _\infty.
\end{eqnarray}

\begin{theorem}
\label{theorem::close}
Suppose that \textup{(A0)--(\ref{eqa2})} hold for $p$
and that \textup{(A0)} holds for~$\tilde p$.
Let $\psi= \max\{ \varepsilon,\varepsilon',\varepsilon''\}$
and let
$\Psi= \max\{ \varepsilon,\varepsilon',\varepsilon'',\varepsilon'''\}$.
When $\Psi$ is sufficiently small:

(1) Conditions \textup{(A1)} and \textup{(\ref{eqa2})} hold for $\tilde p$.

(2) We have:
$\operatorname{\mathsf{Haus}}(R,\tilde R) \leq\frac{2C \psi}{\beta}$.

(3) If $\operatorname{\mathsf{reach}}_\mu(R)>0$ for some $\mu>0$, then
$\tilde{R}\oplus\frac{4\psi}{\mu^2}  \VHF R$.
\end{theorem}

\begin{pf}
(1) Write the spectral decompositions
$H = U \Lambda U^T$ and
$\tilde H = \tilde U \tilde\Lambda\tilde U^T$.
By (\ref{eq::weyl2}),
$|\lambda_j - \tilde\lambda_j| \leq
D \Vert H - \tilde{H}\Vert _{\mathrm{max}} \leq D \varepsilon''$.
Thus, $\tilde p$ satisfies (A1) when $\varepsilon''$ is small enough.
Clearly, (\ref{eqa2}) also holds as long as $\Psi$ is small enough.

(2) By the Davis--Kahan theorem (\ref{eq::davis-kahan}),
\[
\Vert L - \tilde{L}\Vert \leq\frac{\Vert H -
\tilde H\Vert _F}{\beta} \leq\frac{D \Vert H - \tilde H\Vert
 _{\mathrm{max}}}{\beta} \leq
\frac{D
\varepsilon
''}{\beta}.
\]
For each $x$,
\begin{eqnarray*}
\bigl\Vert G(x) - \tilde G(x)\bigr\Vert &=&
\bigl\Vert L(x) g(x) - \tilde{L}(x) \tilde g(x)\bigr\Vert
\\
&\leq& \bigl\Vert \bigl(L(x) - \tilde{L}(x) \bigr) g(x) \bigr\Vert +
\bigl\Vert \tilde{L}(x) \bigl(
\tilde{g}(x) -g(x) \bigr) \bigr\Vert
\\
& \leq& \frac{D \Vert g(x)\Vert  \varepsilon''}{\beta} + \varepsilon'.
\end{eqnarray*}
It follows that,
$\Vert L-\tilde L\Vert \leq C \varepsilon''$ and
$\sup_x \Vert G(x) - \tilde G(x)\Vert  \leq C \psi$.

Now let $\tilde x\in\tilde R$.
Thus, $\Vert \tilde G(\tilde x)\Vert =0$, and hence
$\Vert G(\tilde x)\Vert  \leq C \psi$.
Let $\gamma$ be the path through $\tilde x$ so that
$\gamma(s) = \tilde x$ for some $s$.
Let $r=\gamma(0)\in R$.
From part 2 of Lemma~\ref{lemma::thiss},
note that
$\xi'(s) = \Vert G(\tilde x)\Vert $.
We have
\[
C \psi\ge\bigl\Vert G(\tilde x)\bigr\Vert = \xi'(s) = \xi'(0)
+ s \xi''(u)
\]
for some $u$ between $0$ and $s$.
Since $\xi'(0)=0$,
from part 4 of Lemma~\ref{lemma::thiss},
$C \psi\ge s \xi''(u)\ge\frac{s\beta}{2}$
and so
$\Vert r-\tilde x\Vert  \leq s \leq\frac{2C\psi}{\beta}$.
Thus, $d(\tilde x,R)\leq\Vert r-\tilde x\Vert  \leq2C \psi/\beta$.

Now let $x\in R$.
The same argument shows that
$d(x,\tilde R)\leq2C \psi/\beta$
since (A1) and (\ref{eqa2}) hold for $\tilde p$.

(3) Choose any fixed $\kappa>0$ such that
$\kappa< \frac{\mu^2}{5\mu^2 + 12}$.
When $\Psi$ is sufficiently small,
$\Psi\leq\kappa\operatorname{\mathsf{reach}}_\mu(K)$.
Then
$\tilde R\oplus\frac{4\psi}{\mu^2}  \VHF R$
from Theorem~\ref{thm::CCSL}.
\end{pf}

\section{Ridges of density estimators}

Now we consider estimating the ridges in the density ridge case
(no hidden manifold).
Let $X_1,\ldots, X_n\sim P$
where $P$ has density $p$ and let
%
\begin{equation}
\hat p_h(x) = \frac{1}{n}\sum_{i=1}^n
\frac{1}{h^D} K \biggl(\frac
{\Vert x-X_i\Vert }{h} \biggr)
\end{equation}
be a kernel density estimator
with kernel $K$ and bandwidth $h$.
Let $\hat R$ be the ridge defined by $\hat p$.
In this section, we bound
$\operatorname{\mathsf{Haus}}(R,\hat R)$.
We assume that $P$ is supported on a compact set $\cK\subset\R^D$
and that $p$ and its first, second and third derivatives
vanish on the boundary of $\cK$.
(This ensures there is no boundary bias in the kernel density estimator.)

We assume that all derivatives of $p$ up to and including
fifth degree are bounded and continuous.
We also assume the conditions on the kernel
in Gine and Guillou (\citeyear{gine2002rates})
which are satisfied by all the usual kernels.
Results on
$\Vert p(x) - \hat p_h(x)\Vert _\infty$
are given, for example, in
\citet{prakasa1983nonparametric}, \citet{gine2002rates} and \citet
{yukich1985laws}.
The results in those references imply that
\[
\varepsilon\equiv\sup_{x\in\cK} \bigl\Vert p(x) - \hat p(x)\bigr\Vert
_\infty= O\bigl(h^2\bigr) + O_P \biggl(\sqrt{
\frac{\log n}{n h^D}} \biggr).
\]
For the derivatives,
rates are proved in the sense of mean squared error by
\citet{chacon2011asymptotics}.
They can be proved in the $L_\infty$ norm
using the same techniques as in
\citet{prakasa1983nonparametric}, \citet{gine2002rates} and \citet
{yukich1985laws}.
The rates are:
\begin{eqnarray*}
\varepsilon' &\equiv& \max_j \sup
_{x\in\cK} \bigl|g_j(x) - \hat g_j(x)\bigr| = O
\bigl(h^2\bigr) + O_P \biggl(\sqrt{\frac{\log n}{n h^{D+2}}}
\biggr),
\\
\varepsilon'' &\equiv& \max_{j,k}
\sup_{x\in\cK} \bigl|H_{j,k}(x) - \hat H_{j,k}(x)\bigr| =
O\bigl(h^2\bigr) + O_P \biggl(\sqrt{\frac{\log n}{n h^{D+4}}}
\biggr),
\\
\varepsilon''' &\equiv& \sup
_{x\in\cK} \bigl\Vert H'(x) - \hat H'(x)
\bigr\Vert _{\mathrm{max}} = O\bigl(h^2\bigr) + O_P \biggl(
\sqrt{\frac{\log n}{n h^{D+6}}} \biggr).
\end{eqnarray*}
[See \citet{arias2013estimation}, e.g.]
Let
$\psi_n =  (\frac{\log n}{n} )^{{2}/{(D+8)}}$.
Choosing
$h \asymp\sqrt{\psi_n}$
we get that
$\varepsilon\asymp\varepsilon' \asymp\varepsilon'' \asymp
O_P(\psi_n)   \mbox{ and }  \varepsilon''' = o_P(1)$.
From
Theorem~\ref{theorem::close}
and the rates above we have the following.

\begin{theorem}
\label{thm::main1}
Let $\hat R^* = \hat R \cap(R\oplus\delta)$.
Under the assumptions above
and assuming that \textup{(A1)} and \textup{(\ref{eqa2})} hold,
we have, with
$h \asymp\sqrt{\psi_n}$
that
%
\begin{equation}
\operatorname{\mathsf{Haus}}\bigl(R,\hat R^*\bigr) = O_P(\psi_n).
\end{equation}
If $\operatorname{\mathsf{reach}}_\mu(R)>0$ then
$\hat R^*\oplus O(\psi_n)  \VHF R$.
\end{theorem}

Let $\overline{p}_h(x) = \mathbb{E}(\hat p_h(x))$
and let $R_h$ be the ridge set of $\overline{p}_h$.
It may suffice for practical purposes to
estimate $R_h$ for some small $h>0$.
Indeed, as a corollary to Theorem~9 in the next section
(letting $R$ take the role of $M$ and $R_h$ take the role
of~$R_\sigma$) it follows that
$\Hausdorff(R,R_h) = O(h^2)$ and
$R_h$ is topologically similar to $R$.
In this case, we can take $h$ fixed rather than letting it tend to 0.
For fixed $h$, we then have dimension-independent rates.\looseness=1

\begin{theorem}
\label{thm::main2}
Let $h>0$ be fixed and let
$\tilde\psi_n = \sqrt{\log n/n}$.
Let $\hat R^* = \hat R \cap(R\oplus\delta)$.
Under the assumptions above
and assuming that \textup{(A1)} and \textup{(\ref{eqa2})} hold for $R_h$
we have,
that
%
\begin{equation}
\operatorname{\mathsf{Haus}}\bigl(R_h,\hat R^*\bigr) = O_P(\tilde
\psi_n).
\end{equation}
If $\operatorname{\mathsf{reach}}_\mu(R_h)>0$ then
$\hat R^*\oplus O(\tilde\psi_n)  \VHF R$.
\end{theorem}

\section{Ridges as surrogates for hidden manifolds}
\label{sec::hidden}

Consider now the case where
$P_\sigma= (1-\eta)\operatorname{Unif}(\cK) + \eta(W\star\Phi_\sigma)$
where $W$ is supported on $M$.
We assume in this section that $M$ is a compact manifold without boundary.
We also assume that $W$ has a twice-differentiable density $w$ respect
to the uniform measure on $M$.
(Here, $w$ is a function on a smooth manifold and
the derivatives are defined in the usual way, that is, with respect
to any coordinate chart.)
We also assume that $w$ is bounded away from zero and $\infty$.
In this section, we add the subscript $\sigma$ to
the density, the gradient, etc. to emphasize the dependence on $\sigma$.
For example, the density of $P_\sigma$ is denoted by
$p_\sigma$,
the gradient by $g_\sigma$ and the Hessian by $H_\sigma$.

We want to show that the ridge of $p_\sigma$
is a surrogate for $M$.
Specifically, we show that, as $\sigma$ gets small,
there is a subset $R_*\subset R$ in a neighborhood of $M$ such that
$\operatorname{\mathsf{Haus}}(M,R_*) = O(\sigma^2\log(1/\sigma))$ and such that
$R_* \VHF M$.
We assume that $\eta=1$ in what follows; the extension to
$0<\eta<1$ is straightforward.
We also assume that $M$ is a compact $d$-manifold
with positive reach $\kappa$.
We need to assume that $M$ has positive reach rather than just
positive $\mu$-reach.
The reason is that, when $M$ has positive reach,
the measure $W$ induces a smooth distribution
on the tangent space $T_x M$
for each $x\in M$.
We need this property in our proofs but
this property is lost
if $M$ only has positive $\mu$-reach for some $\mu<1$
due to the presence of unsmooth features such as corners.

The density of $X$ is
%
\begin{equation}
p_\sigma(x) = \int_M \phi_\sigma(x-z)
\,dW(z),
\end{equation}
where
$\phi_\sigma(u) = (2\pi)^{-D/2} \sigma^{-D}$
$\exp (- \frac{\Vert u\Vert ^2}{2\sigma^2} )$.
Thus, $p_\sigma$ is a mixture of Gaussians.
However, it is a rather unusual mixture;
it is a \emph{singular mixture of Gaussians} since the mixing
distribution $W$
is supported on a lower dimensional manifold.

Let $T_x M$ be the tangent space\setcounter{footnote}{3}\footnote{
Recall that the tangent space at a point $x$ is the
linear space spanned by the derivative vectors of smooth curves
on the manifold through that point.}
to $M$ at $x$ and let
$T_x^\perp M$ be the normal space to $M$ at $x$.
Define the \emph{fiber} at $x\in M$ by
$F_x=T_x^\perp M \cap B_D(x,r)$.
A~consequence of the fact that the reach $\kappa$ is positive
and $M$ has no boundary is that,
for any $0 < r < \kappa$,
$M\oplus r$ can be written as a disjoint union
%
\begin{equation}
\label{eq::fibers} M\oplus r = \bigcup_{x\in M}
F_x.
\end{equation}

Let $r_\sigma>0$ satisfy the following conditions:
%
\begin{equation}
r_\sigma< \sigma, \frac{r_\sigma}{\sigma^2}\to\infty, \qquad r_\sigma\log
\biggl(\frac{1}{\sigma^{2+D}} \biggr) = o(1)\qquad \mbox{as } \sigma\to0.
\end{equation}
Specifically, take
$r_\sigma= \alpha\sigma$ for some
$0 < \alpha< 1$.
Fix any $A\ge2$ and define
%
\begin{equation}
K_\sigma= \sqrt{ 2\sigma^2 \log \biggl(\frac{1}{\sigma
^{A+D}}
\biggr)}.
\end{equation}

\begin{theorem}[(Surrogate theorem)]
\label{theorem::main}
Suppose that $\kappa=\operatorname{\mathsf{reach}}(M)>0$.
Let $R_\sigma$ be the ridge set of $p_\sigma$.
Let
$M_\sigma= M \oplus r_\sigma$
and
$R^*_\sigma= R_\sigma\cap M_\sigma$.
For all small $\sigma>0$:
\begin{longlist}[1.]
\item[1.]$R^*_\sigma$ satisfies \textup{(A1)} and \textup{(\ref{eqa2})}
with $\beta= c \sigma^{-(D-d+2)}$ form some $c>0$.
\item[2.]$\operatorname{\mathsf{Haus}}(M,R^*_\sigma) = O(K_\sigma^2)$.
\item[3.]$R^*_\sigma\oplus C K_\sigma^2 \VHF M$.
\end{longlist}
If $R_\sigma$ is instead taken to be the ridge set of $\log p_\sigma$
then the same results are true
with
$\beta= c \sigma^{-2}$ and
$M_\sigma= M \oplus\kappa$.
\end{theorem}

\begin{remark*} Without the assumption that $M$ has no boundary,
there would be boundary effects of order $K_\sigma$.
That is, the Hausdorff distance behaves like
$O(K_\sigma)$ for points near the boundary and like
$O(K_\sigma^2)$ for points not near the boundary.
\end{remark*}

The theorem shows that
in a neighborhood of the manifold,
there is a well-defined ridge,
that the ridge is close to the manifold
and is nearly homotopic to the manifold.
It is interesting to compare the above result to recent work
on finite mixtures of Gaussians
[\citet{carreira2003number,EdFaRosocg2012}].
In those papers, it is shown that there can be fewer or more
modes than the number of Gaussian components in a finite mixture.
However, for small $\sigma$, it is easy to see that for each
component of the mixture, there is a nearby mode.
Moreover, the density will be highly curved at those modes.
Theorem~\ref{theorem::main} can be thought of as a
version of the latter two facts for the case of manifold mixtures.

The theorem refers to the ridges defined by $p_\sigma$ and
the ridges defined by $\log p_\sigma$.
Although the location of the ridge sets is the same for both cases,
the behavior of the function around the ridges is different.
There are several reasons we might want to use $\log p$ rather than $p$.
First, when $p$ is Gaussian,
the ridges of $\log p$ correspond to the usual principal components.
Second, the surrogate theorem holds in an $O(1)$ neighborhood of $M$
for the log-density whereas it only holds in an
$O(\sigma)$ neighborhood of $M$ for the density.

To prove the theorem, we need a preliminary result.
Let
%
\begin{equation}
\tilde\sigma= \sigma\log^3 \biggl(\frac{1}{\sigma^{D+A}} \biggr).
\end{equation}
Given a point $x$ let $\hat x$ be its projection onto $M$.
In what follows,
if $T$ is a matrix, then an expression of the form
$T + O(r_n)$
is to be interpreted to mean
$T+B_n$ where
$B_n$ is a matrix whose entries are of order $O(r_n)$.
Let
%
\begin{equation}
\phi_\perp(u) = \frac{e^{-\Vert u\Vert ^2/(2\sigma^2)}}{(2\pi)^{
{(D-d)}/{2}}\sigma^{D-d}},\qquad u\in\R^{D-d}.
\end{equation}

\begin{lemma}\label{lemma::important-quantities}
For all $x\in M_\sigma$,
\begin{longlist}[1.]
\item[1.]$p_\sigma(x) = \phi_\perp(x-\hat x) (1+ O(\tilde\sigma
) )$.
\item[2.]
Let $p_{\sigma,B}(x) = \int_{M\cap B} \phi_\sigma(x-z)\,dW(Z)$.
Then
$p_{\sigma,B}(x) =
\phi_\perp(x-\hat x) (1+ O(\tilde\sigma) )$.
\item[3.]
$g_\sigma(x) = -\frac{1}{\sigma^2}p_\sigma(x) ( (x-\hat x) +
O(K_\sigma^2)  )$ and
$\Vert g_\sigma(x)\Vert  = O(\sigma^{-(D-d-1)})$.
\item[4.]
The eigenvalues of $H_\sigma(x)$ are
%
\begin{equation}
\lambda_j(x) = %
\cases{ O(\tilde\sigma), &\quad $j \leq d,$
\vspace*{2pt}
\cr
\displaystyle -\frac{p_\sigma(x)}{\sigma^2} \biggl[1 - \frac{d_M^2(x)}{\sigma^2} + O(\tilde
\sigma) \biggr], &\quad $j = d+1,$\vspace*{2pt}
\cr
\displaystyle-\frac{p_\sigma(x)}{\sigma^2} \bigl[1 + O(\tilde
\sigma) \bigr], &\quad $j > d+1.$} %
\end{equation}
%
%
\item[5.] The projection matrix $L_\sigma$ satisfies
\[
L_\sigma(x) = \left[ \pmatrix{ 0_{d} &
0_{d,D-d}
\vspace*{2pt}\cr
0_{D-d,d} & I_{D-d} }
 \right] + O(\tilde
\sigma).
\]
\item[6.] Projected gradient:
\[
G_\sigma(x)= -\frac{1}{\sigma^2} \bigl( (x-\hat x)\phi_\perp(x-
\hat x) \bigl(1+ O(\tilde\sigma)\bigr) + O_\perp \bigl(K^2_\sigma
\bigr) \bigr),
\]
where $O_\perp(K^2_\sigma)$ is a term of size $O(K_\sigma^2)$ in
$T_x^\perp$.
\item[7.] Gap:
\[
\lambda_d(x)-\lambda_{d+1}(x) \ge \frac{p_\sigma(x)}{\sigma^2}
\bigl[1-\alpha^2 + O(\tilde\sigma ) \bigr]
\]
and
\[
\beta\equiv\inf_{x\in R\oplus\delta}\bigl[\lambda_d(x)-
\lambda_{d+1}(x)\bigr] \geq c \sigma^{-(D=d+2)}
\]
and $\lambda_{d+1}(x) \leq-\beta$.
\item[8.]$\Vert H'_\sigma\Vert _{\mathrm{max}}=
O(\sigma^{-(D+3-d)})$.\vadjust{\goodbreak}
\end{longlist}
\end{lemma}

\begin{pf}
The proof is
quite long and technical and so we relegate
it to the \hyperref[app]{Appendix}.
\end{pf}

\begin{pf*}{Proof of Theorem~\ref{theorem::main}}
Let us begin with the ridge based on $p_\sigma$.

(1) Condition (A1) follows
from parts 8 and 1 of Lemma~\ref{lemma::important-quantities}
together with equation
(\ref{eq::perp}).

To verify (\ref{eqa2}), we use
parts 3 and 8 of Lemma~\ref{lemma::important-quantities}:
we get, for all small $\sigma$, that
\begin{eqnarray*}
\Vert L^\perp g\Vert\bigl \Vert H'\bigr\Vert _{\mathrm{max}} &\leq&
\Vert g\Vert \bigl\Vert H'\bigr\Vert _{\mathrm{max}} \leq \frac{c}{\sigma^{D-d-1}}
\frac{1}{\sigma^{D+3-d}} < \frac{c^2\sigma^2 }{\sigma^{2(D-d+2)}} \leq \sigma^2
\beta^2 \\
&\leq& \frac{\beta^2}{2 D^{3/2}}
\end{eqnarray*}
as required.

(2) Suppose that $x\in R_\sigma^*$.
Then $\Vert G_\sigma(x)\Vert =0$.
Let $\hat x$ be the unique projection of $x$ onto $M$.
From part 6 of Lemma~\ref{lemma::important-quantities},
\[
\bigl\Vert (x-\hat x)\phi_\perp(x-\hat x) \bigl(1+O(\tilde\sigma)\bigr) +
O_\perp \bigl(K^2\bigr)\bigr\Vert =0,
\]
and hence $x = \hat x + O(K_\sigma^2)$.

Now let $\hat x \in M$.
From the expression above, we see that
$\Vert G_\sigma(\hat x)\Vert =O(K_\sigma^2)$.
Let $\gamma$ be the path through $x$ and let $r$ be the destination of
the path.
Hence $\gamma(s) =x$ for some $s$ and $\gamma(0) = r$.
Now we use Lemma~\ref{lemma::thiss}.
Then $\Vert G\Vert  = \xi'$ and
\[
O\bigl(K_\sigma^2\bigr) = \xi'(s) =
\xi'(s) - \xi'(0) = s\xi''(
\tilde s) \ge\Vert x-\hat r\Vert \xi''(\tilde s) \ge
\Vert x-\hat r\Vert \beta/2
\]
and so
$\Vert x-\hat r\Vert  = O(K_\sigma^2)$.
Hence, $\operatorname{\mathsf{Haus}}(R_\sigma,M) = O(K_\sigma^2)$.

(3) Homotopy.
This follows from part (2)
and Theorem~\ref{thm::CCSL}.

Now consider the ridges of
$\log p_\sigma(x)$.
The proof is essentially the same as the proof above.
The main difference is the Hessian as we now explain.
Note that the Hessian $H_\sigma^*$ for $\log p_\sigma(x)$ is
\[
H_\sigma^*(x) = \frac{1}{p_\sigma(x)} \biggl( H_\sigma(x) -
\frac{1}{p_\sigma(x)}g_\sigma(x) g_\sigma ^T(x)
\biggr).
\]
From
Lemma~\ref{lemma::important-quantities}, parts 3 and 4,
it follows that (after an appropriate rotation),
\[
H_\sigma^*(x) = -\frac{1}{\sigma^2} \biggl( \left[ \matrix{
 O_{d} & 0_{d \times D-d}
\vspace*{2pt}\cr
0_{D-d \times d} & I_{D-d} }
 \right] + O(\tilde
\sigma) \biggr).
\]
Hence,
\[
\lambda_{d+1}(x) = -\frac{1}{\sigma^2} + O(\tilde\sigma)
\]
and
\[
\lambda_{d+1}(x) -\lambda_d(x) = \frac{1}{\sigma^2} + O(
\tilde \sigma).
\]
Notice in particular, that the dominant term of the smallest eigenvalue of
$-\beta H_\sigma^*(x)$ is 1 whereas
that the dominant term of the smallest eigenvalue of
$-\beta H_\sigma(x)$ is 1 $d^2_M(x)/\sigma^2$
which is why we required $\Vert x-\hat x\Vert $ to be less than $\sigma$ in
Theorem~\ref{theorem::main}.
Here, we only require that
$\Vert x-\hat x\Vert  \leq\kappa$.
\end{pf*}

We may now combine
Theorems \ref{theorem::close},
\ref{thm::main1},
\ref{thm::main2} and
\ref{theorem::main} to get the following.

\begin{corollary}
Let $\hat R^*$ be defined as in
Theorem~\ref{thm::main1}.
Then
%
\begin{equation}
\operatorname{\mathsf{Haus}}\bigl(\hat R^*,M\bigr) = O_P \biggl( \biggl(
\frac{\log n}{n} \biggr)^{{2}/{(D+8)}} \biggr) + O\bigl(K_\sigma^2
\bigr).
\end{equation}
Similarly, if $\hat R^*$ be defined as in
Theorem~\ref{thm::main2}
then
%
\begin{equation}
\operatorname{\mathsf{Haus}}\bigl(\hat R^*,M\bigr) = O_P \biggl(\sqrt{
\frac{\log n}{n}} \biggr) + O\bigl(K_\sigma^2
+h^2\bigr).
\end{equation}
\end{corollary}

\section{SuRFing the ridge}
\label{sec::SCMS}

Here, we discuss
Subspace Ridge Finding (SuRF)
by using
density estimation,
followed by denoising
and then followed by
the subspace constrained mean shift (SCMS) algorithm
due to \citet{Principal}.
We will not go into great details about the algorthm; we refer the
reader to
\citet{Principal}.

Let us begin by reviewing the mean shift algorithm.
The \emph{mean shift algorithm}
[\citet{Fukunaga,cheng1995mean,Comaniciu}]
is a method for finding the modes of a density
by approximating the steepest ascent paths.
The algorithm starts with a mesh of points and then moves the points along
the gradient ascent trajectories toward local maxima.

Given a sample $X_1,\ldots, X_n$ from $p$,
consider the kernel density estimator
%
\begin{equation}
\hat{p}_h(x) = \frac{1}{n}\sum_{i=1}^n
\frac{1}{h^D} K \biggl(\frac{\Vert x-X_i\Vert }{h} \biggr),
\end{equation}
where $K$ is a kernel and $h>0$ is a bandwidth.
Let
$\cM= \{{v_1,\ldots, v_m}\}$ be a collection of mesh points.
These are often taken to be the same as the data
but in general they need not be.
Let $v_j(1)=v_j$ and
for $t=1,2,3,\ldots$ we define the trajectory
$v_j(1),v_j(2),\ldots,$ by
%
\begin{equation}
v_j(t+1) = \frac{\sum_{i=1}^n X_i K
( {\Vert v_j(t) - X_i\Vert }/{h} )}{
\sum_{i=1}^n K ( {\Vert v_j(t) - X_i\Vert }/{h} )}.
\end{equation}

It can be shown that each trajectory
$\{v_j(t)\dvtx t=1,2,3,\ldots,\}$
follows the gradient ascent path and
converges to a mode of $\hat{p}_h$.
Conversely,
if the mesh $\cM$ is rich enough, then
for each mode of $\hat{p}_h$, some trajectory will converge to that mode.

\begin{figure}
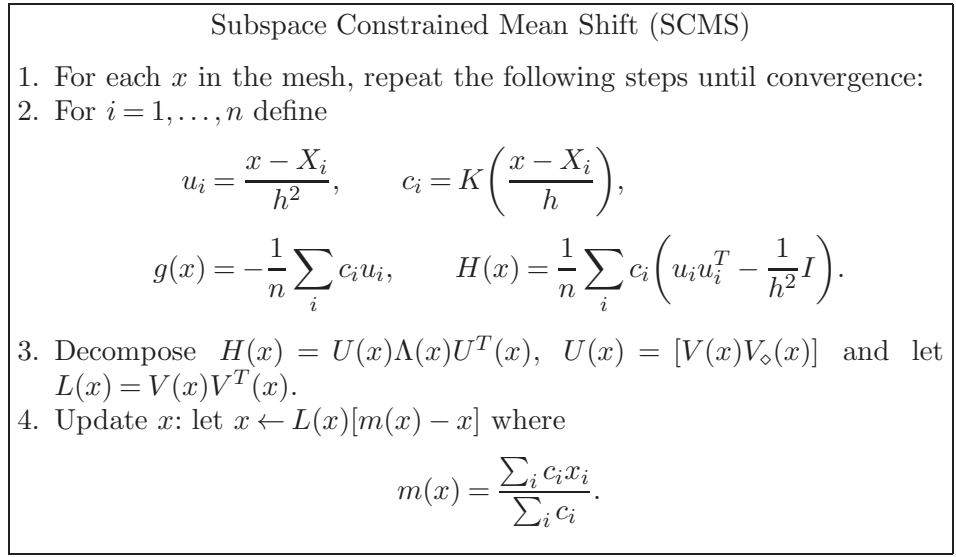

\fbox{\parbox{4.85in}{
\begin{center}
Subspace Constrained Mean Shift (SCMS)
\end{center}
\begin{enumerate}
\item For each $x$ in the mesh, repeat the following steps until convergence:
\item For $i=1,\ldots, n$ define
\begin{eqnarray*}
u_i &=& \frac{x-X_i}{h^2},\qquad c_i = K \biggl(
\frac{x-X_i}{h} \biggr),
\\
g(x) &=& -\frac{1}{n}\sum_i
c_i u_i,\qquad H(x) = \frac{1}{n}\sum
_i c_i \biggl( u_i
u_i^T - \frac{1}{h^2}I \biggr).
\end{eqnarray*}
\item Decompose $H(x) = U(x) \Lambda(x) U^T(x)$,
$U(x) = [V(x)  V_\diamond(x)]$ and let $L(x) = V(x) V^T(x)$.
\item Update $x$: let
$x \leftarrow L(x) [m(x) -x]$
where
\[
m(x)= \frac{\sum_i c_i x_i}{\sum_i c_i}.
\]
\end{enumerate}
}}
\caption{SCMS algorithm from Ozertem and Erdogmus (\citeyear{Principal}).}
\label{fig::scmsalg}
\end{figure}

\begin{figure}[b]
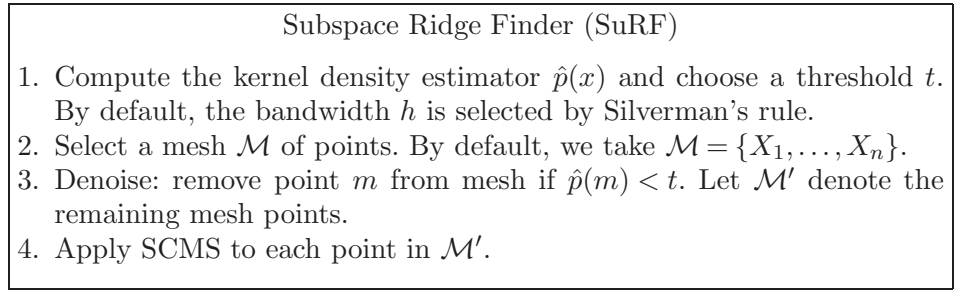

\fbox{\parbox{4.85in}{
\begin{center}
Subspace Ridge Finder (SuRF)
\end{center}
\begin{enumerate}
\item Compute the kernel density estimator $\hat p(x)$ and choose a
threshold $t$.
By default, the bandwidth $h$ is selected by Silverman's rule.
\item Select a mesh ${\cal M}$ of points.
By default, we take ${\cal M} = \{X_1,\ldots, X_n\}$.
\item Denoise: remove point $m$ from mesh if $\hat p(m) < t$.
Let ${\cal M}'$ denote the remaining mesh points.
\item Apply SCMS to each point in ${\cal M}'$.
\end{enumerate}
}}
\caption{SuRF algorithm.}
\label{fig::surfalg}
\end{figure}

The SCMS algorithm mimics the mean shift algorithm
but it replaces the gradient with the projected
gradient at each step.
The algorithm can be applied to $\hat p$ or any
monotone function of $\hat p$.
As we explained earlier, there are some advantages to using $\log\hat p$.
Figure~\ref{fig::scmsalg}
gives the algorithm
for the log-density.
This is the version we will use
in our examples.
Figure~\ref{fig::surfalg}
gives the full SuRF algorithm.

The SCMS algorithm provides a numerical approximation to
the paths $\gamma$ defined by the projected gradient.
We illustrate the numerical algorithm in Section~\ref{section::examples}.

\section{Implementation and examples}
\label{section::examples}

Here, we demonstrate ridge estimation
in some two-dimensional examples.
In each case, we will find the one-dimensional ridge set.
Our purpose is to show proof of concept;
there are many interesting implementation details that we will not
address here.
In each case, we use SuRF.

\begin{figure}

\includegraphics{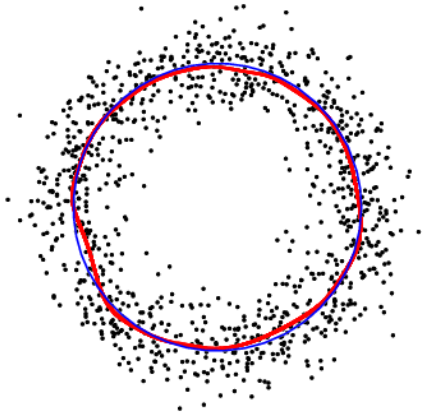}

\caption{Estimated hyper-ridge set (red curve) from data generated
from a circular manifold $M$ (blue curve) of radius 3.
The sample size is 1000, using Normal noise with $\sigma= 0.5$.
The estimate is computed from a kernel density estimator using
the Silverman Normal reference rule for the bandwidth.
The starting points for the modified SCMS algorithm are taken
the evaluation points of the density estimator excluding the
points below 25\% of the maximum estimated density.}
\label{fig::good-circle}
\end{figure}

To implement the method requires that we choose a bandwidth $h$
for the kernel density estimator.
There has been recent work on bandwidth selection
for multivariate density estimators
such as
\citeauthor{chacon2010multivariate}
(\citeyear{chacon2010multivariate,chacon2012bandwidth}) and
\citet{panaretos2012nonparametric}.
For the purposes of this paper, we simply use the
Silverman rule [\citet{scott1992multivariate}].

\begin{figure}[t!]

\includegraphics{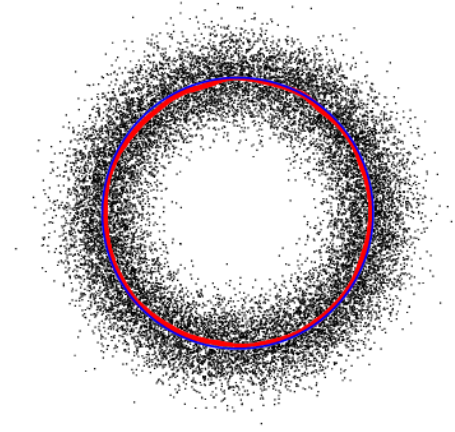}

\caption{Estimated hyper-ridge set (red curve) from data generated
from a circular manifold $M$ (blue curve) of radius 3.
The sample size is 20,000, using Normal noise with $\sigma= 0.5$.
The estimate is computed from a kernel density estimator using
the Silverman Normal reference rule for the bandwidth.
The starting points for the modified SCMS algorithm are taken
the evaluation points of the density estimator excluding the
points below 25\% of the maximum estimated density.}
\label{fig::good-circle-thb}\vspace*{-12pt}
\end{figure}

\begin{figure}[t!]

\includegraphics{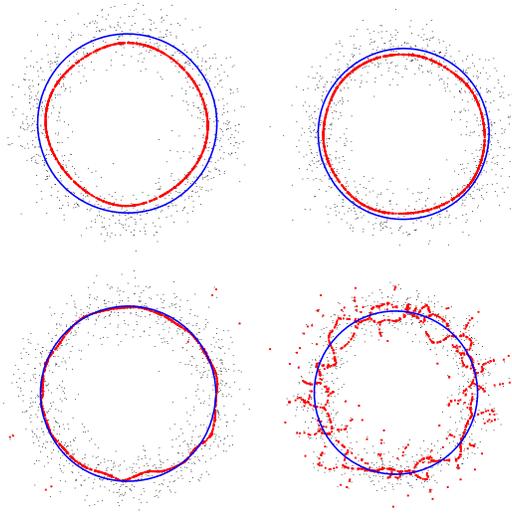}

\caption{Effect of decreasing bandwidth.
The data are i.i.d. samples from the same manifold as in the previous figure.
Eventually we reach a phase transition
where the structure of the estimator falls apart.}
\label{fig::good-circles-phase}
\end{figure}

\begin{figure}

\includegraphics{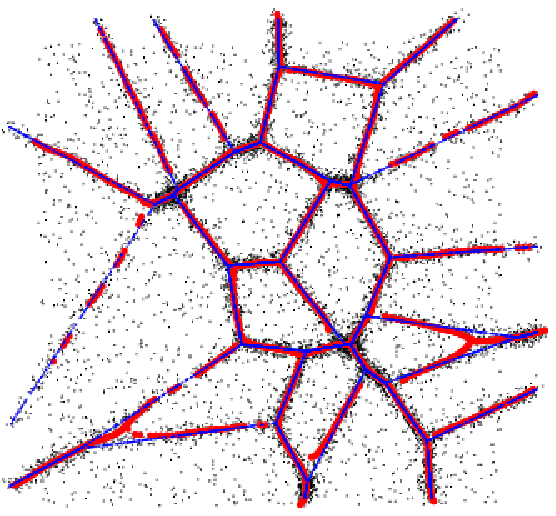}

\caption{Data generated from a stylized ``cosmic web'' consisting of
intersecting
line segments and a uniform background clutter.
Total sample size is 10,000. The starting points for the modified SCMS
algorithm are taken
the evaluation points of the density estimator excluding the
points below 5\% of the maximum estimated density.}
\label{fig::good-web-th}
\end{figure}

Figures~\ref{fig::good-circle} through \ref{fig::good-web-th} show two examples of SuRF.
In the first example, the manifold is a circle.
Although the circle example
may seem easy,
we remind the reader that no existing
statistical algorithms that we are aware of
can, without prior assumptions,
take a point cloud as input and find a circle,
automatically.

The second example
is a stylized ``cosmic web''
of intersecting line segments and with random background clutter.
This is a difficult\ case that violates the assumptions; specifically
the underlying object does not have positive reach.
The starting points for the SCMS algorithm are a subset of the grid points
at which a kernel density estimator is evaluated. We select those points
for which the estimated density is above a threshold relative to the
maximum value.

Figure~\ref{fig::good-circles-phase} shows
the estimator for four bandwidths.
This shows an interesting phenomenon.
When the bandwidth $h$ is large, the estimator is
biased (as expected) but it is still
homotopy equivalent to the true $M$.
However, when $h$ gets too small,
we see a phase transition where the estimator falls
apart and degenerates into small pieces.
This suggests it is safer to oversmooth and have a small
amount of bias.
The dangers of undersmoothing are greater
than the dangers of oversmoothing.

The theory in Section~\ref{sec::hidden}
required the underlying structure to have positive reach
which rules out intersections and corners.
To see how the method fares when these assumptions are violated,
see Figure~\ref{fig::good-web-th}.
While the estimator is far from perfect,
given the complexity of the example,
the procedure does surprisingly well.

\section{Conclusion}

We presented an analysis of nonparametric ridge
estimation.
Our analysis had two main components:
conditions that guarantee that the estimated ridge
converges to the true ridge, and
conditions to relate the ridge to an underlying hidden manifold.

We are currently investigating several questions.
First, we are finding the minimax rate for this problem
to establish whether or not our proposed method is optimal.
Also, \citet{klemela2005adaptive}
has derived mode estimation procedures that adapt to the local
regularity of the mode.
It would be interesting to derive similar adaptive theory for ridges.
Second, the hidden manifold case required that the manifold had
positive reach.
We are working on relaxing this condition to allow for
corners and intersections
(often known as stratified spaces).
Third, we are developing an extension where
ridges of each dimension $d=0,1,\ldots$ are found sequentially
and removed one at a time.
This leads to a decomposition of the point cloud into
structures of increasing dimension.
Finally, there are a number of methods for speeding up
the mean shift algorithm.
We are investigating how to adapt these speedups for
SuRF.

As we mentioned in the \hyperref[sec1]{Introduction},
there is recent work on metric graph reconstruction
which is a way of modeling intersecting filaments
[\citet{aanjaneya2012metric,lecci2013statistical}].
These algorithms have the advantage of being designed to handle
intersecting ridges.
However, it appears that they are very sensitive to noise.
Currently, we are investigating the idea of first
running SuRF and then applying metric graph reconstruction.
Preliminary results suggest that this approach may get
the best of both approaches.

\begin{appendix}
\section*{Appendix}\label{app}

The purpose
of this appendix
is to prove
Lemma~\ref{lemma::important-quantities}.
Recall that
the gradient is
$g_\sigma(x) = - \frac{1}{\sigma^2}\int_M (x-z)\phi_\sigma(x-z) \,dW(z)$
and the Hessian is
%
\begin{eqnarray}
H_\sigma(x) &=&
-\frac{1}{\sigma^2} \int_M
\biggl( I - \frac{ (x-z)(x-z)^T}{\sigma^2} \biggr)\phi _\sigma(x-z) \,dW(z)
\nonumber
\\[-8pt]
\\[-8pt]
\nonumber
&=& -\frac{1}{\sigma^2} \biggl[ p_\sigma(x) I - \frac{1}{\sigma^2}\int
_M (x-z) (x-z)^T\phi _\sigma (x-z)
\,dW(z) \biggr].
\end{eqnarray}

We can partition $M_\sigma$
into disjoint fibers.
Choose an $x\in M_\sigma$ and let
$\hat x$ be the unique projection of $x$ onto $M$.
Let $B = B(\hat x,K_\sigma)$.
For any bounded function $f(x,z)$,
%
\begin{equation}
\label{eq::tail}\quad  \int_{M\cap B^c}f(x,z) \phi_\sigma(x-z)
\,dW(z) \le \frac{C}{(2\pi)^{{D}/{2}}}\frac{e^{-K^2_\sigma/(2\sigma
^2)}}{\sigma
^{D}} W\bigl(B^c\bigr) \leq
C\sigma^{A}.
\end{equation}

Let $T=T_{\hat x} M$ denote the $d$-dimensional
tangent space at $\hat x$ and let
$T^\perp$ denote the
$(D-d)$-dimensional normal space.
For $z\in B\cap M$,
let $\overline{z}$ be the projection of $z$ onto $T$.
Then
%
\begin{equation}
\label{eq::decompose} x -z = (x-\hat x) + (\hat x - \overline{z}) + (\overline{z}-z) =
d_M(x)u + (\hat x - \overline{z}) + R,
\end{equation}
where
$u = (x-\hat x)/d_M(x) \in T^\perp$
and
$R = (\overline{z}-z)$.
[Recall that $d_M$ is the distance function; see (\ref{eq::distfun}).]
For small enough $\sigma$,
there is a smooth map $h$ taking $z$ to $\overline{z}$ that is a bijection
$B\cap M$ and so
the distribution $W$
induces a distribution $\overline{W}$, that is,
$\overline{W}(A) = W(h^{-1}(A))$.
Let $\overline{w}$ denote the density of $\overline{W}$
with respect to Lebesgue measure $\mu_d$ on $T$.
The density is bounded above and below
and has two continuous derivatives.

\begin{lemma}
For every $x\in R\oplus\sigma$,
$\sup_{z\in B}\Vert z-\overline{z}\Vert  \leq c K^2_\sigma$.
\end{lemma}

\begin{pf}
Choose any $z\in B$ and let
$\overline{z}$ be its projection onto $T$.
Because the reach is $\kappa>0$,
there exists a ball
$S(a,\kappa)\subset\R^D$
such that
$a$ is in the plane defined $\hat x, z$ and $\overline{z}$,
$S(a,\kappa)$ is tangent to the manifold at $\hat x$
and
$S(a,\kappa)$
does not intersect $M$ except at $\hat x$.
Consider the line through
$z$ and $\overline{z}$ and let $z_a$
be the point where the line intersects $S(a,\kappa)$.
Now
$\Vert z_a-\overline{z}\Vert \ge\Vert z-\overline{z}\Vert $
and by elementary geometry,
$\Vert z_a-\overline{z}\Vert \leq C K^2_\sigma$.
\end{pf}

Recall that $r_\sigma= \alpha\sigma$ with $0 < \alpha< 1$.
Define the following quantities:
\begin{eqnarray*}
\beta&=& \frac{e^{-\alpha^2/2}(1-\alpha^2)}{2\sigma^{D-d+2}},
\qquad \tilde\sigma= \sigma\log^3 \biggl(
\frac{1}{\sigma^{D+A}} \biggr),
\\
\phi_\perp(u) &=& \frac{e^{-\Vert u\Vert ^2/(2\sigma^2)}}{(2\pi)^{
{(D-d)}/{2}}\sigma^{D-d}},\qquad \phi_\Vert (w) =
\frac{e^{-\Vert w\Vert ^2/(2\sigma^2)}}{(2\pi)^{{d}/{2}}\sigma^{d}},
\\
p_{\sigma,B}(x) &=& \int_{M\cap B} \phi_\sigma(x-z)
\,dW(z),\qquad B = B(\hat x,K_\sigma),
\end{eqnarray*}
where $u\in\R^{D-d}$ and $w\in\R^d$.

\begin{lemma}
We have that
%
\begin{equation}
\phi_\sigma(x-z) = \phi_\perp(x-\hat x)\phi_\Vert (
\hat x - \overline {z}) \bigl(1+ O(\tilde\sigma)\bigr).
\end{equation}
\end{lemma}

\begin{pf}
First note that, for all $x\in R_\sigma$,
%
\begin{equation}
\label{eq::perp} \frac{1}{(2\pi)^{ {(D-d)}/{2} }} \frac{e^{-\alpha^2/2}}{\sigma^{D-d}} \leq \phi_\perp(x-
\hat x) \leq \frac{1}{(2\pi)^{{(D-d)}/{2}}}\frac{1}{\sigma^{D-d}}
\end{equation}
and so,
$\phi_\perp(x-\hat x) \asymp\sigma^{-(D-d)}$
as $\sigma\to0$.
Now,
\[
\Vert x-z\Vert ^2 = \Vert x-\hat x\Vert ^2 + \Vert \hat
x - z\Vert ^2 + \Vert \overline z - z\Vert ^2 + 2\langle
x-\hat x,\overline z - z\rangle,
\]
we have that
\[
\phi_\sigma(x-z) = \phi_\perp(x-\hat x) \phi_\Vert (
\hat x-\overline z) e^{-\Vert z-\overline{z}\Vert ^2/(2\sigma^2)} e^{ -\langle x-\hat x,\overline z - z\rangle/\sigma^2}.
\]
Now
$\Vert z-\overline{z}\Vert ^2=O(K_\sigma^4)$ and
$|\langle x-\hat x,\overline z - z\rangle| \leq
\Vert  x-\hat x\Vert  \Vert \overline z - z\Vert  = O(\sigma K^2_\sigma)$
and so
\[
e^{-\Vert z-\overline{z}\Vert ^2/(2\sigma^2)} e^{ -\langle x-\hat x,\overline z - z\rangle/\sigma^2} = \bigl(1+ O(\tilde\sigma)\bigr).
\]
\upqed\end{pf}

\begin{pf*}{Proof of Lemma~\ref{lemma::important-quantities}}
{1}. From (\ref{eq::tail}),
$p_\sigma(x) = \int_{M\cap B} \phi_\sigma(x-z) \,dW(z) + O(\sigma^A)$.
Now
\begin{eqnarray*}
&&\int_{M\cap B} \phi_\sigma(x-z) \,dW(z) \\
&&\qquad= \bigl(1+ O(
\tilde\sigma)\bigr) \phi_\perp(x-\hat x)\int_{M\cap B}
\phi_\Vert (\hat x - \overline {z}) \,dW(z)
\\
&&\qquad= \bigl(1+ O(\tilde\sigma)\bigr) \phi_\perp(x-\hat x)\int
_{T M\cap B} \phi_\Vert (\hat x - \overline {z})
\overline{w}(\overline{z}) \,d\mu_d(\overline{z})
\\
&&\qquad= \bigl(1+ O(\tilde\sigma)\bigr) \phi_\perp(x-\hat x) \int
_{T} \frac{1}{(2\pi)^{d/2}}e^{-\Vert t\Vert ^2/2} \overline{w}(\hat x
+ \sigma t) \,d\mu_d(t),
\end{eqnarray*}
where
$A= \{ t=(\overline{z}-\hat x)/\sigma\dvtx  \overline{z}\in B\}$.
The volume of
$T$ is $O(\sigma^{D+A})$
and $T\to\R^d$ as $\sigma\to0$.
Also, $\overline{w}(\hat x + \sigma t) = \overline{w}(\hat x) +
O(\sigma)$.
Hence,
\[
\int_{T} \overline{w}(\hat x + \sigma t) \,d
\mu_d(t) = \bigl(\overline{w}(\hat x)+O(\sigma)\bigr) \bigl( 1 - O
\bigl(\sigma^{D+A}\bigr)\bigr)
\]
and so
\[
\int_{M\cap B} \phi_\sigma(x-z) \,dW(z) =
\phi_\perp(x-\hat x) \bigl(1+ O(\tilde\sigma)\bigr) \bigl(\overline{w}(
\hat x)+O(\sigma)\bigr) \bigl( 1 - O\bigl(\sigma^{D+A}\bigr)\bigr)
\]
and
\begin{eqnarray*}
p_\sigma(x) &=& \phi_\perp(x-\hat x) \bigl(1+ O(\tilde\sigma)
\bigr) \bigl(\overline{w}(\hat x)+O(\sigma)\bigr) \bigl( 1 - O\bigl(
\sigma^{D+A}\bigr)\bigr) + O\bigl(\sigma ^A\bigr)
\\
&=& \phi_\perp(x-\hat x) \bigl(1+ O(\tilde\sigma)\bigr).
\end{eqnarray*}

{2}. \emph{$p_{\sigma,B}(x)$.} This follows since in part 1 we
showed that
$p_{\sigma,B}(x) = p_{\sigma}(x) +O(\sigma^A)$.

{3}. For the gradient, we have
\begin{eqnarray*}
-\sigma^2 g_\sigma(x) &= &\int(x-z) \phi_\sigma(x-z)
\,dW(z)
\\
&=& (x-\hat x)\int\phi_\sigma(x-z) \,dW(z)\\
&&{} + \int(\hat x - \overline {z})
\phi_\sigma(x-z) \,dW(z) \\
&&{}+ \int(\overline{z} - z)\phi_\sigma(x-z)
\,dW(z)
\\
&=& \mathrm{I} + \mathrm{II} + \mathrm{III}.
\end{eqnarray*}
Now,
$\mathrm{I} = (x-\hat x)p_\sigma(x) = (x-\hat x)\phi_\perp(x-\hat x) (1+
O(\tilde\sigma))$
and
\begin{eqnarray*}
\mathrm{II} &=& \int_{M\cap B} (\hat x - \overline{z})
\phi_\sigma(x-z) \,dW(z) + O\bigl(\sigma ^A\bigr)
\\
&=& \bigl(1+ O(\tilde\sigma)\bigr)\phi_\perp(x-\hat x)\int
_{M\cap B} (\hat x - \overline{z})\phi_\Vert (\hat x-
\overline{z}) \,dW(z) + O\bigl(\sigma^A\bigr).
\end{eqnarray*}
For some $u$ between $\hat x$ and $\overline{z,}$ we have
\begin{eqnarray*}
&&\int_{M\cap B} (\hat x - \overline{z})\phi_\Vert (\hat
x-\overline {z}) \,dW(z)\\
&&\qquad= \sigma\int_{M\cap B} \frac{\hat x - \overline{z}}{\sigma}
\phi _\Vert (\hat x-\overline{z}) \,dW(z)
\\
&&\qquad= \sigma\int_{h^{-1}(B)} \frac{\hat x - \overline{z}}{\sigma} \phi_\Vert
(\hat x- \overline{z})\overline{w}(\overline{z}) \,d\mu _d(
\overline{z})
\\
&&\qquad= \sigma\int_{A} t \frac{1}{(2\pi)^{d/2}}e^{-\Vert t\Vert ^2/2}
\overline {w}(\hat x + \sigma t) \,d\mu_d(t)
\\
&&\qquad= \sigma\int_{A} t \frac{1}{(2\pi)^{d/2}}e^{-\Vert t\Vert ^2/2}
\bigl[\overline{w}(\hat x)+ \overline{w}'(\hat x)\sigma t +
\overline{w}''(u)\sigma^2
t^2/2\bigr] \,d\mu_d(t)
\\
&&\qquad= O\bigl(\sigma^2\bigr),
\end{eqnarray*}
where
$A = \{ t = (\overline{z}-\hat x)/\sigma\in h^{-1}(B)\}$.
Finally,
\begin{eqnarray*}
\mathrm{III}& =& \int_{M\cap B} (\overline{z} - z)
\phi_\sigma(x-z) \,dW(z) + O\bigl(\sigma ^A\bigr)\\
& =&
\phi_\perp(x-\hat x)O\bigl(K_\sigma^2\bigr) + O
\bigl(\sigma^A\bigr)\\
& =& O\bigl(K_\sigma ^2\bigr)
\phi _\perp(x-\hat x).
\end{eqnarray*}
Hence,
\begin{eqnarray*}
-\sigma^2 g_\sigma(x) &=& (x-\hat x)p_\sigma(x) + O
\bigl(\sigma^2\bigr) + \phi_\perp(x-\hat x)O
\bigl(K_\sigma^2\bigr) \\
&=& p_\sigma(x) \bigl((x-\hat
x) + O\bigl(K_\sigma^2\bigr)\bigr)
\end{eqnarray*}
and hence
\[
g_\sigma(x) = -\frac{1}{\sigma^2}p_\sigma(x) \bigl( (x-\hat x)
+ O\bigl(K_\sigma^2\bigr) \bigr).
\]
It follow from part 1 that $\Vert g_\sigma(x)\Vert  = O(\sigma^{-(D-d-1)})$.

{4}.
To find the eigenvalues, we first approximate the
Hessian.
Without loss of generality,
we can rotate the coordinates so that
$T$ is spanned by $e_1,\ldots, e_d$,
$T^\perp$ is spanned by $e_{d+1},\ldots, e_D$
and $u =(0,\ldots, 0,1)$.
Now,
\[
-\frac{\sigma^2 H_\sigma(x)}{p_\sigma(x)} = I - \frac{ \int(x-z)(x-z)^T \phi_\sigma(x-z) \,dW(z)}{\sigma^2 \int
\phi
_\sigma(x-z) \,dW(z)}
\]
and
\[
\int_{M\cap B^c} (x-z) (x-z)^T \phi_\sigma(x-z)
\,dW(z) = O\bigl(\sigma^A\bigr).
\]
Let $Q = \int_{M\cap B} (x-z)(x-z)^T \phi_\sigma(x-z) \,dW(z)$.
Then, from (\ref{eq::decompose}), we have
$Q = Q_1 + Q_2 + Q_3 + Q_4 + Q_5 + Q_6$
where
\begin{eqnarray*}
Q_1 &=& d_M^2(x)u u^T\int
_{M\cap B} \phi_{\sigma}(x-z) \,dW(z),\\ Q_2& =& \int
_{M\cap B} (\hat x - \overline{z}) (\hat x - \overline{z})^T
\phi_\sigma(x-z) \,dW(z),
\\
Q_3 &=& \int_{M\cap B} (\overline{z}-z) (
\overline{z}-z)^T \phi _\sigma (x-z) \,dW(z),\\ Q_4 &=&
\int_{M\cap B} (x - \hat x) (\hat x-\overline{z})^T
\phi _\sigma (x-z) \,dW(z),
\\
Q_5 &= &\int_{M\cap B} (x - \hat x) (
\overline{z}-z)^T \phi_\sigma(x-z) \,dW(z),\\ Q_6 &=&
\int_{M\cap B} (\hat x-\overline{z}) (\overline{z}-z)^T
\phi _\sigma(x-z) \,dW(z).
\end{eqnarray*}

First, we note that
\[
Q_1 = d_M^2(x)u u^T
\phi_\perp(x-\hat x) \bigl(1+ O(\tilde\sigma)\bigr).
\]
Next,
\begin{eqnarray*}
Q_2 &=& \int_{M\cap B} (\hat x - \overline{z}) (\hat
x - \overline{z})^T \phi_\sigma(x-z) \,dW(z)
\\
&= &\bigl(1+ O(\tilde\sigma)\bigr)\phi_\perp(\hat x-x) \int
_{M\cap B} (\hat x - \overline{z}) (\hat x - \overline{z})^T
\phi _\Vert (\hat x - \overline{z}) \,dW(z)
\end{eqnarray*}
and
\begin{eqnarray*}
&&\int_{M\cap B} (\hat x - \overline{z}) (\hat x -
\overline{z})^T \phi _\Vert (\hat x - \overline{z}) \,dW(z) \\
&&\qquad=
\int_{h^{-1}(B)} (\hat x - \overline{z}) (\hat x -
\overline{z})^T \phi _\Vert (\hat x - \overline{z})
\overline{w}(\overline{z}) \,d\mu_d(\overline{z})
\\
&&\qquad= \overline{w}(\hat x) \int_{h^{-1}(B)} (\hat x - \overline{z}) (
\hat x - \overline{z})^T \phi _\Vert (\hat x - \overline{z})
\,d\mu_d(\overline{z})+O\bigl(K_\sigma^5\bigr).
\end{eqnarray*}
Next,
with $t=(t_1,\ldots,t_d,0,\ldots, 0)$,
\begin{eqnarray*}
&&\int_{h^{-1}(B)} (\hat x - \overline{z}) (\hat x -
\overline{z})^T \phi _\Vert (\hat x - \overline{z}) \,d
\mu_d(\overline{z})\\
&&\qquad = \sigma^2 \int_{\overline{B}}
t t^T (2\pi)^{-d/2}e^{-\Vert t\Vert ^2/2} \,d\mu _d(t)
\\
&&\qquad= \sigma^2 \biggl( \int t t^T (2\pi)^{-d/2}e^{-\Vert t\Vert ^2/2}
\,d\mu_d(t) \\
&&\hspace*{16pt}\qquad\quad{}- \int_{\overline
{B}^c} t t^T (2
\pi)^{-d/2}e^{-\Vert t\Vert ^2/2} \,d\mu_d(t) \biggr)
\\
&&\qquad= \sigma^2 \biggl(\left[ \matrix{
I_d & 0
\vspace*{2pt}\cr
0 & 0 }
 \right] + O\bigl(\sigma^{A+D}\bigr)
\biggr)
\end{eqnarray*}
and so
\begin{eqnarray*}
Q_2 &=& \bigl(1+O(\tilde\sigma)\bigr)\phi_\perp(x-\hat x)
\sigma^2\\
&&{}\times \biggl(\left[ \matrix{ I_d &
0
\vspace*{2pt}\cr
0 & 0 }
 \right] + O\bigl(\sigma^{A+D}\bigr)
\biggr).
\end{eqnarray*}

A similar analysis on the remaining terms yields:
\begin{eqnarray*}
Q_3 &=& \bigl(1+O(\tilde\sigma)\bigr)\phi_\perp(x-\hat x) O
\bigl(K_\sigma^4\bigr),
\\
Q_4 &=& \bigl(1+O(\tilde\sigma)\bigr)\phi_\perp(x-\hat x) O
\bigl(\sigma K_\sigma ^2\bigr),
\\
Q_5 &=& \bigl(1+O(\tilde\sigma)\bigr)\phi_\perp(x-\hat x) O
\bigl(\sigma K_\sigma ^2\bigr),
\\
Q_6 &= &\bigl(1+O(\tilde\sigma)\bigr)\phi_\perp(x-\hat x) O
\bigl(K_\sigma^3\bigr).
\end{eqnarray*}

Combining all the terms, we have
\begin{eqnarray*}
Q &=& \bigl(1+ O(\tilde\sigma)\bigr)\phi_\perp(x-\hat x)\\
&&{}\times \left(
d_M^2(x)u u^T + \sigma^2 \biggl(
\left[ %
\matrix{I_d & 0
\vspace*{2pt}\cr
0 & 0 }
 \right] + O\bigl(\sigma^{A+D}\bigr)
\biggr) \right) \\
&&{}+ O \bigl(\sigma K^2_\sigma \bigr).
\end{eqnarray*}
Hence,
\begin{eqnarray*}
H_\sigma(x) &=& - \bigl(1+O(\tilde\sigma)\bigr)\frac{p_\sigma(x)}{\sigma^2}
\\
&&{}\times \left( \left[ %
\matrix{ 0 & \cdots& 0 & 0 &
\cdots& \cdots& \cdots& 0
\vspace*{2pt}\cr
\vdots& \ddots& \vdots& 0 & \cdots& \cdots& \cdots& 0
\vspace*{2pt}\cr
0 & \cdots& 0 & 0 & \cdots& \cdots& \cdots& 0
\vspace*{2pt}\cr
\hline 0 & \cdots& 0 & 1 & 0 & \cdots& \cdots& 0
\vspace*{2pt}\cr
0 & \cdots& 0 & 0 & 1 & \cdots& \cdots& 0
\vspace*{2pt}\cr
0 & \cdots& 0 & 0 & 0 & \ddots& \cdots& 0
\vspace*{2pt}\cr
0 & \cdots& 0 & 0 & 0 & \cdots& 1 & 0
\vspace*{2pt}\cr
0 & \cdots& 0 & 0 & 0 & \cdots& 0 & 1-\frac{d_M^2(x)}{\sigma^2}
} %
 \right] +O (\tilde\sigma ) \right).
\end{eqnarray*}
The result follows.

{5}. This follows from part 4 and the Davis--Kahan theorem.

{6}. From part 5,
$L_\sigma(x) = L^\dagger+ E$ where
$L^\dagger=
[
{0_{d\times d} \enskip 0_{d,D-d}\atop
0_{D-d,d} \enskip I_{D-d}
}
]$
and
$E = O (\tilde\sigma )$.
Hence,
$G_\sigma(x) = L_\sigma(x) g_\sigma(x) =(L^\dagger+ E) g_\sigma(x)$
and the result follows from parts~3 and 4.

{7}. These follow from part 4.

{8}. Now we turn to $\Vert H'_\sigma\Vert $.
Let $\Delta= (x-z)$.
We claim that
\begin{eqnarray*}
H' &=& \frac{1}{\sigma^4}\int \biggl[(\Phi\otimes I) ( I\otimes
\Delta + \Delta\otimes I) \\
&&\hspace*{30pt}{}- \frac{\phi_\sigma(\Delta)}{\sigma^2}\bigl(I\otimes\Delta
\Delta^T\bigr) \bigl( \operatorname{\mathsf{vec}}(I) \otimes\Delta^T\bigr)
\biggr]\,dW(z)
\\
&&{} + \frac{1}{\sigma^4} \int\phi_\sigma(\Delta) \bigl( \operatorname{\mathsf{vec}}(I)
\otimes \Delta^T\bigr) \,dW(z).
\end{eqnarray*}
To see this, note first that
$H= \frac{1}{\sigma^4}Q - \frac{1}{\sigma^2}A$
where
%
\begin{equation}
Q = \int(x-z) (x-z)^T \phi_\sigma(x-z) \,dW(z) \quad\mbox{and}\quad A =
p_\sigma(x) I.
\end{equation}
Note that
$Q = \int(x-z)(x-z)^T \Phi \,dW(z)$
where $\Phi= \phi_\sigma(\Delta) I_D$.
So
\[
Q' = \int(d/dx)\bigl[(x-z) (x-z)^T \Phi\bigr] \,dW(z)
\]
and
\begin{eqnarray*}
\frac{d}{dx} \bigl[(x-z) (x-z)^T \Phi\bigr]& =&
\frac{d (x-z)(x-z)^T \Phi}{dx}\\
& =& \frac{ d \Delta\Delta^T \Phi}{d\Delta}.
\end{eqnarray*}
Now
$(d/dx)( \Delta\Delta^T \Phi) = (fg)'$
where
$f = \Delta\Delta^T$ and $g = \Phi$ and so
\begin{eqnarray*}
\frac{d}{dx}\bigl( \Delta\Delta^T \Phi\bigr) &=& (\Phi\otimes
I) \frac{d}{dx}\bigl( \Delta\Delta^T\bigr) + \bigl(I\otimes
\Delta \Delta ^T\bigr) \frac{d}{dx} \Phi
\\
&=& (\Phi\otimes I) ( I\otimes\Delta + \Delta\otimes I) - \frac{\phi_\sigma(\Delta)}{\sigma^2}
\bigl(I\otimes\Delta\Delta^T\bigr) \bigl( \operatorname{\mathsf{vec}}(I) \otimes
\Delta^T\bigr).
\end{eqnarray*}
Hence,
\[
Q' = \int \biggl[(\Phi\otimes I) ( I\otimes\Delta + \Delta\otimes
I) - \frac{\phi_\sigma(\Delta)}{\sigma^2}\bigl(I\otimes\Delta\Delta^T\bigr) \bigl( \operatorname{\mathsf{vec}}(I) \otimes\Delta^T\bigr) \biggr]\,dW(z).
\]
By a similar calculation,
\[
A' = -\frac{1}{\sigma^2} \int\phi_\sigma(\Delta) \bigl( {
\sf vec}(I) \otimes\Delta^T\bigr) \,dW(z).
\]
Thus,
\begin{eqnarray*}
H' &=& \frac{1}{\sigma^4} Q' - \frac{1}{\sigma^2}
A'
\\
&=& \frac{1}{\sigma^4}\int \biggl[(\Phi\otimes I) ( I\otimes\Delta + \Delta
\otimes I) \\
&&\hspace*{32pt}{}- \frac{\phi_\sigma(\Delta)}{\sigma^2}\bigl(I\otimes\Delta\Delta^T\bigr)
\bigl( \operatorname{\mathsf{vec}}(I) \otimes\Delta^T\bigr) \biggr]\,dW(z)
\\
&&{} + \frac{1}{\sigma^4} \int\phi_\sigma(\Delta) \bigl( \operatorname{\mathsf{vec}}(I)
\otimes \Delta^T\bigr) \,dW(z).
\end{eqnarray*}
Each of these terms is of order
$O(\sup_{x\in M}\Vert w''(x)\Vert /\sigma^{D-d+1})$.
Consider the first term
\begin{eqnarray*}
\frac{1}{\sigma^4}\int(\Phi\otimes I) ( I\otimes\Delta) \,dW(z) &=&
\frac{1}{\sigma^{4+D}} (2\pi)^{D/2} \int e^{-\Vert x-z\Vert ^2/(2\sigma^2)} ( I\otimes\Delta)
\,dW(z)
\\
&=& \frac{1}{\sigma^{3+D} (2\pi)^{D/2} }\int e^{-\Vert u\Vert ^2/2} ( I\otimes u) \,dW(z),
\end{eqnarray*}
where $u = (x-z)/\sigma$.
As in the proof of part 1,
we can restrict to $B\cap M$, do a change of measure to $\overline W$
and the term is dominated by
\begin{eqnarray*}
&&\frac{1}{\sigma^{3+D-d}(2\pi)^{D/2} } \int_A e^{-\Vert u\Vert ^2/2} ( I\otimes u)
\bigl[\overline{w}(\hat x) + \overline{w}'(\tilde u)\sigma u \bigr]
\,d\mu_d(t)
\\
&&\qquad=\frac{C}{\sigma^{3+D-d} (2\pi)^{D/2} }.
\end{eqnarray*}
\upqed\end{pf*}
The other terms may be bounded similarly.
\end{appendix}

\section*{Acknowledgements}
The authors thank the reviewers for
many suggestions that improved the paper.
In particular, we thank the Associate Editor
who suggested a simplified proof of Lemma~\ref{lemma::thiss}.

%



\printaddresses
\end{document}